\documentclass[a4paper,11pt]{amsart}

\usepackage[a4paper, left=2.5cm, right=3cm, top=2.5cm, bottom=3cm]{geometry}
\usepackage{amsmath,amssymb,amsthm,amsfonts,mathtools}
\usepackage{chngcntr}
\usepackage[american]{babel}
\counterwithin*{equation}{section}
\counterwithin*{equation}{section}

\usepackage[hidelinks]{hyperref}
\usepackage{todonotes}
\usepackage{cite}
\theoremstyle{plain} \numberwithin{equation}{section}
\newtheorem{theorem}{Theorem}[section]
\newtheorem{proposition}[theorem]{Proposition}
\newtheorem{lemma}[theorem]{Lemma}
\newtheorem{corollary}[theorem]{Corollary}
\theoremstyle{definition}
\newtheorem{definition}[theorem]{Definition}
\newtheorem{remark}[theorem]{Remark}

\title{Vector Valued Katok-Sarnak Formulas}
\author{Rebekka Strathausen}
\address{Fachbereich Mathematik, Technische Universität Darmstadt, Schlossgartenstrasse 7,
64289 Darmstadt, Germany}
\email{strathausen@mathematik.tu-darmstadt.de}

\begin{document}
\maketitle
\begin{abstract} We introduce a version of the Katok-Sarnak formula linking twisted traces of scalar valued Maass cusp forms for $\Gamma_0(N)$ with $N$ square-free to Fourier coefficients of vector valued Maass cusp forms for the metaplectic group $\textup{Mp}_2(\mathbb{Z})$ and the Weil representation associated to a lattice of signature $(2,1)$. We prove this by computing the Fourier coefficients of the vector valued Siegel and Millson theta lift of a scalar valued newform by two different methods and comparing the results. \end{abstract}

\section{Introduction}
\subsection{Motivation}
A well known theorem by Duke, \cite{MR931205}, states that Heegner points of discriminant $D$ equidistribute in the modular curve $Y := \textup{SL}_2(\mathbb{Z})\setminus\mathbb{H}$ for $|D| \rightarrow \infty$.  Denote by $\Lambda_D$ the finite set of Heegner points corresponding to a discriminant $D < 0$. We furthermore denote the number of such points, i.e. $\textup{SL}_2(\mathbb{Z})$-equivalence classes of $\mathcal{Q}_D^+$, the set of positive definite quadratic forms of discriminant $D$, by $h(D)$. The proof of said theorem follows from Weyl's criterion: For any spectral component $f \in L^2(Y)$, it must be shown that the expression
\begin{align*}
P_D(f) := \frac{1}{h(D)}\sum\limits_{z \in \Lambda_D}f(z)
\end{align*}
goes to $0$ for $|D| \rightarrow \infty$. We call $P_D(f)$ the $D$-th period of $f$. The spectral components in question are given by either Eisenstein series or Maass cusp forms. The periods of Maass cusp forms are related to Fourier coefficients of half integral Maass forms by Katok-Sarnak formulas (see, for example, Theorem \ref{theins} below), which makes it possible to study periods of Maass cusp forms and their asymptotic behaviour for large discriminants in great detail. The original Katok-Sarnak formula was introduced by Katok and Sarnak in their 1993 paper \cite{MR1244668}. It  has since been generalized to periods twisted by a genus character and to Maass forms of even weight for $\Gamma_0(N)$ with $N$ odd in \cite{MR3549627} and \cite{MR4389022}. It can be applied to prove refinements of Duke's Theorem, as done in \cite{MR3549627}. \\
\indent The goal of the present article is to prove a version of the Katok-Sarnak formula which describes twisted periods of Maass cusp forms of even weight for the congruence subgroup $\Gamma_0(N)$ for square-free $N$ in terms of Fourier coefficients of vector valued Maass cusp forms of half-integral weight for the Weil representation associated to a lattice of signature $(2,1)$ and the metaplectic group $\textup{Mp}_2(\mathbb{Z})$.\\
%%%%%%%%%%%%%%%%%%%%%%%%%%%%%%%%%%%%%%%%%%%%%%%%%%%%%
\indent Both of the aforementioned generalizations, \cite{MR3549627} and \cite{MR4389022}, of the Katok-Sarnak formula use Shimura's original construction of the Shimura correspondence and do not apply methods of theta lifting as Katok and Sarnak did. The aim of this paper is to reprove and extend the current results by means of theta lifting, considering a vector valued instead of a scalar valued setting in order to limit computational expenses and obtain a more structural understanding of the correspondence between integral and half-integral weight Maass cusp forms.\\
\indent This approach has the following advantages: First, it makes it possible to state a Katok-Sarnak formula for all square-free levels, not limited to the level being odd. Second, it is an approach easily generalized to different settings, since vector valued theta lifts can be defined for other contexts as well, such as, for example, for Hilbert or Siegel modular forms.

\subsection{Maass Cusp Forms}
A \textit{Maass cusp form of weight $k$} for a congruence subgroup $\Gamma$ of $\textup{SL}_2(\mathbb{Z})$ is a real analytic  function $f : \mathbb{H} \rightarrow \mathbb{C}$, which transforms like a modular form of weight $k$ for  $\Gamma$ and is an eigenfunction of the weight $k$ Laplace operator
\begin{displaymath}
\Delta_k := -y^2\left(\frac{\partial^2}{\partial x^2} + \frac{\partial^2}{\partial y^2}\right) + kiy\left(\frac{\partial}{\partial x} + i\frac{\partial}{\partial y} \right),
\end{displaymath}
where $z = x + iy$.
Additionally, we request that $f(z)$ is of exponential decay in the cusps of $\Gamma$.\\
\indent Denote the eigenvalue of a Maass cusp form under $\Delta_k$ by $\lambda = (s - k/2) (1-s-k/2)$ with $s \in \mathbb{C}$. For fixed  $s$, the Maass cusp forms of weight $k$ with eigenvalue $\lambda$ form a finite dimensional vector space. Any element $f$ of this space has a Fourier expansion of the form
\begin{displaymath}
%\sum\limits_{h \in L'/L}\big[\sum\limits_{\substack {n \in \mathbb{Z} + q(h) \\n \neq 0}}c(h, n)(4\pi |n|v)^{-\frac{k}{2}}W_{\textup{sign}(n)\frac{k}{2}, s- \frac{1}{2}}(4\pi |n|v)e(nx)\big]\mathfrak{e}_h
f(z) = \sum_{n \not= 0}a(n)(4\pi |n|y)^{-\frac{k}{2}}W_{\textup{sign}(n)\frac{k}{2}, s-\frac{1}{2}}(4\pi |n|y)e(nx),
\end{displaymath}
where $W_{\mu,\nu}$ is the $W$-Whittaker function. We call a Maass cusp form even if its Fourier coefficients fulfill the relation
\begin{align}\label{ev}
a(-n) = \frac{\Gamma\left(s + \frac{k}{2}\right)}{\Gamma\left(s - \frac{k}{2}\right)}a(n),
\end{align}
where $\Gamma$ is the standard $\Gamma$-function.
%even hinten wieder raus

\subsection{Twisted Traces}

\indent Now, for simplicity, set $\Gamma = \textup{SL}_2(\mathbb{Z})$. To any fundamental discriminant $\Delta$, we assign a genus character $\chi_{\Delta}$ acting on quadratic forms, as follows.\\
\indent For $D \in \mathbb{Z}$ a discriminant, $\mathcal{Q}_D$ is the set of binary integral quadratic forms $Q = [a, b, c]$ with $Q(x,y) = ax^2 +bxy +cy^2$ of discriminant $D = b^2 -4ac$.  We denote by $\overline{\Gamma}_Q$ the stabilizer of a quadratic form $Q$ in $\overline{\Gamma} = \textup{PSL}_2(\mathbb{Z}) = \Gamma\slash\{\pm 1\}$. It is finite if $D < 0$ and trivial or infinite cyclic if $D > 0$, depending on whether $D$ is a square or not.\\
\indent For $ D < 0$, we can associate to $Q \in \mathcal{Q}_D$ a CM point $z_Q \in \mathbb{H}$, defined as the unique point in $\mathbb{H}$ with $Q(z_Q,1)=0$,  and for $D > 0$, we can associate to $Q \in \mathcal{Q}_D$ a geodesic $C_Q$ in $\mathbb{H}$, defined by the solutions of $a|z|^2 + bx + c = 0$.\\
\indent Let now $\Delta$ be a fundamental discriminant. Following \cite{MR0909238}, we define the \textit{genus character} associated to $\Delta$ in the following sense. For a quadratic form $Q = [a, b, c] $, we pick an integer $n$ coprime to $\Delta$ which is represented by $Q$. Then the genus character $\chi_{\Delta}(Q)$ is given by
\begin{align*}
\chi_{\Delta}(Q) :=  \begin{cases}  \begin{pmatrix}\frac{\Delta}{n} \end{pmatrix}, &\textup{ if } \Delta \mid D, \frac{D}{\Delta} \textup{ is a square mod } 4 \textup { and gcd}(a, b, c, \Delta) =1,\\
0, &\textup{ otherwise}. \\
\end{cases}
\end{align*}
\indent For a weight $0 \leq k$ Maass cusp form $\varphi$ and for $D <0$ we define the twisted traces of CM-values of $\varphi$ by
\begin{align*}
Tr_{\Delta}(\varphi, D) = \sum\limits_{Q\in \Gamma\setminus \mathcal{Q}^+_{|\Delta|D}}\chi_{\Delta}(Q)\frac{1}{|\overline{\Gamma}_Q|}L_{k,z}^{k/2}\varphi(z_Q),
\end{align*}
and for $D > 0$, we define the twisted traces of geodesics by
\begin{align*}
Tr_{\Delta}(\varphi, D) = \sum\limits_{Q\in \Gamma\setminus \mathcal{Q}_{|\Delta|D}}\chi_{\Delta}(Q)|\Delta D|^{\frac{1}{2}-\frac{k}{4}}\int\limits_{C_Q}\varphi(z)Q(z,1)^{\frac{k}{2}-1}dz.
\end{align*}

\noindent Here $L_{k,z} = -2iy^2\frac{\partial}{\partial \overline{z}}$ is the so-called Maass lowering operator, a differential operator which preserves modularity in general but lowers the weight by $2$. The expression $L^{k/2}_{k,z}$ is an iterated version, defined by $L_{2, z} \circ ... \circ L_{k-2, z}\circ L_{k,z}$. For $k < 0$, the definition of the twisted traces remains the same for $D>0$, while for $D < 0$ we exchange $L_{k,z}$ with the raising operator $R_{k,z} = 2i\frac{\partial}{\partial z} + ky^{-1}$.\\

\subsection{Shimura Correspondence}
\indent One of the most well-known correspondences between integral and half-integral modular forms is the so-called Shimura correspondence, which generalizes to Maass cusp forms. It can be realized as a theta lift. The Shimura correspondence assigns to any half-integral weight Maass cusp form in the Kohnen plus space a Maass cusp form of integral weight. It is realized in the following way. Let $k \in 2\mathbb{Z}$ and consider a Maass cusp form $F(\tau)$ of weight $\frac{k+1}{2}$ whose $m$-th Fourier coefficient we denote by $b(m)$ with $b(1) \neq 0$. For $n > 0$, we define $a(n)$ by
\begin{align}\label{shm}
    \sum\limits_{m\vert n} m^{k/2-1}b\left( \frac{n^2}{m^2}\right) = a(n)b(1).
\end{align}
Then the function
\begin{align}\label{shmz}
    f(z) = \sum_{n \not= 0}a(n)(4\pi |n|y)^{-\frac{k}{2}}W_{\textup{sign}(n)\frac{k}{2}, s-\frac{1}{2}}(4\pi |n|y)e(nx)
\end{align}
is a Maass cusp form of weight $k$. Here, we have set $a(-n)$ such that \eqref{ev} is fulfilled.\\
\indent This correspondence was first described by Shimura in \cite{MR332663} and, as a theta lift, can be realized using the Siegel and Millson theta functions (see, for example, \cite{MR0364106}, \cite{MR0389772} or \cite{MR716787}).\\

\subsection{Main Results and Method of Proof}
We present a special case of our main result, see Theorem \ref{KSmain}.
\begin{theorem}\label{theins}
Let $\varphi$ be an even, normalized Maass cusp form of weight $0 \leq k$ for $\Gamma$ that is an eigenform for all Hecke operators. Furthermore, let $\{F_j\}$ be a basis of Hecke eigenforms of the Kohnen plus space of weight $\frac{k+1}{2}$ for $\Gamma' := \Gamma_0(4)$ with the Fourier coefficients of $F_j$ being given by $b_j(n)$. Then we have
\begin{align*}
&\sum\limits_{Shim(F_j) = \varphi}(-1)^{\lceil k/4\rceil}2|\Delta D|^{\frac{1}{2}} \overline{b_j\left((-1)^{k/2}\Delta\right)}b_j\big((-1)^{k/2}D\big)\\
&= (\varphi, \varphi)^{-1}|\Delta D|^{\frac{k}{4}}\begin{cases}
2\sqrt{\pi}Tr_{\Delta}\left(\varphi, D \right) & \textup{ for } \Delta D < 0,\\
2^{\frac{k}{2}}Tr_{\Delta}\left(\varphi, D\right) & \textup{ for } \Delta,D > 0.
\end{cases}
\end{align*}
Here, $Shim(F_j)$ denotes the normalized Shimura lift.
\end{theorem}
\begin{remark}
The main result of this paper are two vector valued Katok-Sarnak formulas, stated in Theorem \ref{KSmain} and Theorem \ref{KSsub}, covering all odd weights, square-free levels and signs of $D,$ and $\Delta$. The full versions, including proof, can be found in \autoref{kks}. The statement of \ref{theins} covers the first two of six possible cases of Theorem \ref{KSmain}.
\end{remark}

The proof of Theorem \ref{theins} is done by computing the Fourier coefficients of a weight $\frac{k+1}{2}$ vector valued Maass form $F$ by two different methods and comparing the results. This Maass form is the Siegel or Millson theta lift of $\varphi$, i.e \\
\begin{align}\label{thlif}
F(\tau) = \int\limits_{\Gamma\setminus\mathbb{H}}\varphi(z)\overline{\Theta(z,\tau)}y^k\frac{dxdy}{y^2},
\end{align}
where $\Theta(z,\tau)$ is the Siegel or Millson theta function.
Here, the integral runs over a fundamental domain for the action of $\Gamma$ on $\mathbb{H}$. The theta kernel $\Theta(z, \tau)$ is modular of weight $k$ in $z$ for $\Gamma$ and modular of weight $\ell$ in $\tau$ for some congruence subgroup $\Gamma'$ of $\textup{SL}_2(\mathbb{Z})$. Typically, $k$ is an integer and $\ell$ is half-integral. For $\varphi$ a Maass cusp form, the integral \eqref{thlif} converges and $F(\tau)$ is modular of weight $\ell$ for $\Gamma'$ with respect to the variable $\tau$. The adjoint theta lift is given by the integral in $\tau$ of a weight $\ell$ modular function for $\Gamma'$ against $\Theta(z,\tau)$. The adjoint lifts of the Siegel and Millson lift we call Shimura and Millson-Shimura lift. They result in a multiple of the function in \eqref{shmz}, derived by Shimura correspondence.\\
\indent We first give some brief overview on Maass cusp forms in Section 2. In Section 3, we introduce the Siegel and Millson theta function. Then, in Section 4, we define the Siegel and Millson theta lift and compute its Fourier coefficients, while in Section 5 we compute the Fourier expansion of the adjoint lifts. Finally, in section 6, we state and prove our main Theorems, the vector valued Katok-Sarnak formulas by bringing together the results from the previous sections.

\subsection{Aknowledgements}
The authors thank Claudia Alfes, Jens Funke, Özlem Imamoglu, Riccardo Zuffetti, Nicolas Potthast, Sören Sprehe and Jan Bruinier for their advice, helpful comments and suggestions. \\
The research of the author is supported by the DFG Collaborative Research Centre TRR 326 “Geometry
and Arithmetic of Uniformized Structures”, project number 444845124 and TRR 358 "Integral Structures in Geometry and Representation Theory", project-ID 491392403.

%%%%%%%%%%%%%%%%%%%%%%%%%%%%%

\section{Preliminaries}\label{sig21}
Let $N \in \mathbb{N}_{>0}$. We attribute to $N$ the quadratic space
\begin{align*}
(V, q) := \biggl(\biggl\{ X := \begin{pmatrix}x_2 & x_1\\ x_3 & -x_2 \end{pmatrix}\Bigm| x_1, x_2, x_3 \in \mathbb{Q} \biggr\},  \hspace{1mm} q(X) := -N\textup{det}(X)\biggl).
\end{align*}
The bilinear form on $V$ corresponding to $q$ is $(X,Y) = N\textup{tr}(XY)$. This space is of signature $(2,1)$ (see Proposition \ref{isom} and Lemma \ref{isomm} below). The group $\textup{GL}^+_2(\mathbb{Q})$ acts isometrically on $V$ via $gX = gXg^{-1}$, for $X \in V$ and $g \in \textup{GL}^+_2(\mathbb{Q})$. Within this space, we regard the lattice
\begin{align*}
 L := \biggl\{\begin{pmatrix}b & -a/N\\ c & -b\end{pmatrix} \Bigm| a, b, c \in \mathbb{Z} \biggr\},
\end{align*}
which is an even lattice in $V$ with corresponding dual lattice given by
\begin{align*}
 L' := \biggl\{\begin{pmatrix}b/2N & -a/N\\ c & -b/2N\end{pmatrix} \Bigm| a, b, c \in \mathbb{Z} \biggr\}.
\end{align*}
We have $L'/L \cong \mathbb{Z}/2N\mathbb{Z}$. The group $\Gamma_0(N)$ acts trivially on $L'/L$ and takes $L$ to itself. $L$ has level $4N$. Set $V(\mathbb{R}) := \mathbb{V} \otimes \mathbb{R}$. By $\mathbb{D}$ we denote the Grassmannian  of lines in $V(\mathbb{R})$ on which $q$ is negative definite:
\begin{align*}
\mathbb{D} = \{Z \subset V(\mathbb{R}) \mid \textup{ dim } Z = 1 \textup{ and } q|_Z < 0 \}. %Quelle aus BIF oder BrO
\end{align*}
The set $\mathbb{D}$ can be identified with the complex upper half plane $\mathbb{H}$ in the following sense (see for example \cite{MR1491448}, Section 11): Let $z = x + iy \in \mathbb{H}$. To $i \in \mathbb{H}$, we associate the line generated by
\begin{align*}
X_Z(i) = \frac{1}{\sqrt{2N}}\begin{pmatrix}0 & 1\\ -1 & 0 \end{pmatrix}.
\end{align*}
The stabilizer of $X_Z(i)$ in $\textup{SL}_2(\mathbb{R})$ is given by $K = \textup{SO}(2)$. We choose $g_z \in \textup{SL}_2(\mathbb{R})$ such that $g_zi = z$ and put
\begin{align*}
X_Z = X_Z(z) := g_zX_Z(i) = \frac{1}{\sqrt{2N}y}\begin{pmatrix}-x & |z|^2\\ -1 & x\end{pmatrix} \in V(\mathbb{R}).
\end{align*}
This gives us an  $\textup{SL}_2(\mathbb{R})$-equivariant isomorphism $\mathbb{H} \rightarrow \mathbb{D}, z \mapsto g_zX_Z(i) = \mathbb{R}X_Z(z)$. Furthermore, we have $(X_Z(z), X_Z(z)) = -1$.\\
\begin{proposition}\label{isom}
Define
\begin{align*}
&X_1(z) := \frac{1}{\sqrt{2N}y}\begin{pmatrix}y & -2xy\\ 0 & -y \end{pmatrix},
\hspace{2mm}\textup{ and }\hspace{2mm}
&X_2(z) := \frac{1}{\sqrt{2N}y}\begin{pmatrix}x & -x^2+y^2\\ 1 & -x \end{pmatrix}.
\end{align*}
Then $X_1(z), X_2(z)$ and $X_3 = X_Z(z)$ form an orthogonal basis of $V(\mathbb{R})$ with $(X_2(z), X_2(z)) = (X_1(z), X_1(z)) = 1$ and $(X_3, X_3) = -1$.
\end{proposition}
\begin{proof}
This can be checked by a direct calculation.
\end{proof}
\begin{lemma}\label{isomm}
Define $\sigma: V(\mathbb{R}) \rightarrow \mathbb{R}^{2,1}$ by
\begin{align*}
\sum\limits_{i =1}^{3}\alpha_iX_i(z) \mapsto (\alpha_1, \alpha_2, \alpha_3). %Vektor umbenennen!!!!!!!
\end{align*}
Then $\sigma$ is an isometry and, consequently, we have $ V(\mathbb{R}) \cong \mathbb{R}^{2,1}$ under $\sigma$.
\end{lemma}
\begin{proof}
This follows directly from Proposition \ref{isom}.
\end{proof}

\subsection{Scalar Valued Maass Cusp Forms}
Let $z = x + iy \in \mathbb{H}$ and $\Gamma = \Gamma_0(N)$ for $N\geq 1$. We define the \textbf{hyperbolic Laplace operator} of weight $k$ by
\begin{displaymath}
\Delta_k := -y^2\bigg(\frac{\partial^2}{\partial x^2} + \frac{\partial^2}{\partial y^2}\bigg) + kiy\bigg(\frac{\partial}{\partial x} + i\frac{\partial}{\partial y} \bigg).
\end{displaymath}
For $k \in \mathbb{Z}$, we define the weight $k$ \textbf{slash operator} on a function $f:\mathbb{H} \rightarrow \mathbb{C}$ and for any matrix  $M = \begin{psmallmatrix}a&b\\c&d\end{psmallmatrix} \in \textup{SL}_2(\mathbb{Z})$ by
\begin{displaymath}
f|_kM(z) = (cz+d)^{-k}f(Mz) =  (cz+d)^{-k}f\biggl(\frac{az+b}{cz+d}\biggr).
\end{displaymath}
We also sometimes denote the term $(cz + d)$ by $j(M, z)$.

\begin{definition} \label{kakadu}
A real analytic function $f: \mathbb{H} \rightarrow \mathbb{C}$ is called a \textbf{Maass wave form} of weight $k$ for $\Gamma$ if it satisfies the following conditions:
\begin{itemize}
\item[i.)] $\Delta_kf = \lambda f$ for some $\lambda \in \mathbb{C}$,
\item[ii.)] $f|_kM(z) = f(z)$ for all $M  \in \Gamma$,
\item[iii.)] There exists $C > 0$ such that for every cusp $\ell$ of $\Gamma\setminus\mathbb{H}$ and $\sigma_{\ell}\in \textup{SL}_2(\mathbb{Z})$ such that $\ell = \sigma_{\ell}\infty$, one has $f|_k\sigma_{\ell} = \mathcal{O}(y^C)$ for $y \rightarrow \infty$.%muss das in jeder Spitze?
\end{itemize}
We call a Maass wave form $f(z)$ a \textbf{Maass cusp form} if, additionally, we have:
\begin{itemize}
\item[iv.)] The expression
\begin{align*}
\int\limits_0^1(f|_k\sigma_{\ell})(x+iy)dx,
\end{align*}
is equal to $0$ for all $y>0$ and all cusps $\ell$.
\end{itemize}
The eigenvalue $\lambda$ is usually normalized by $\lambda = (s - k/2) (1-s-k/2)$ with $s \in \mathbb{C}$. We call $s$ the \textbf{spectral point} of $f(z)$.
\end{definition}
From now on, we focus only on Maass cusp forms. For a given $s$, the Maass cusp forms of weight $k$ for $\Gamma$ with corresponding eigenvalue form a finite dimensional vector space, which we denote by $U_{k,N}(s)$
\begin{proposition}\label{bo}
Let $f(z)$ be a Maass cusp form of weight $k \in \mathbb{Z}$ and its eigenvalue be given by $\lambda = (s-k/2)(1-s-k/2)$. Then $f$ has a Fourier expansion of the form
\begin{displaymath}
f(z) = \sum_{n \not= 0}a(n)(4\pi |n|y)^{-\frac{k}{2}}W_{\textup{sign}(n)\frac{k}{2}, s-\frac{1}{2}}(4\pi |n|y)e(nx),
\end{displaymath}
where $e(q) = e^{2\pi i q}$ and $W_{\mu, \nu}(y)$ is the $W$-Whittaker function \textup{(see Chapter 7 of \textup{\cite{MR232968}})}, which decays exponentially for $y \rightarrow \infty$.
\end{proposition}
\begin{remark}
Condition iv.) of Definition \ref{kakadu} is equivalent with the $0$-th Fourier coefficient $a(0)$ vanishing.
\end{remark}
 \begin{definition} Let $p \nmid N$ be prime. Let  $\varphi \in U_{k,N}(s)$ and its $n$-th Fourier coefficient be given by $a(n)$. Then the action of the $p$-th Hecke operator $T_p$ on $\varphi$ is given by
        \begin{align*}
            &T_p\varphi(z) = p^{k-1}\varphi(pz) + \frac{1}{p}\sum\limits_{j=0}^{p-1}\varphi\Big(\frac{z+j}{p}\Big)\\
            &=\sum_{n \neq 0}\big[a(np) + p^{k-1}a(n/p) \big](4\pi |n|y)^{-\frac{k}{2}}W_{\textup{sign}(n)\frac{k}{2}, s-\frac{1}{2}}(4\pi |n|y)e(nx),
        \end{align*}
        where $a(n/p) = 0$ if $p \nmid n$. If $q|N$ is prime, we define the $q$-th Hecke operator $U_q$ by
        \begin{align*}
            U_q\varphi(z) = \frac{1}{q}\sum\limits_{j=0}^{q-1}\varphi\Big(\frac{z+j}{q}\Big).
        \end{align*}
    \end{definition}

Hecke operators are self-adjoint with respect to the Petersson inner product and commute with $\Delta_k$. There exists an orthogonal basis of $U_{k,N}(s)$ which consists of simultaneous eigenforms for all $T_n$ with $(n,N)=1$. As for modular forms, there is a newform theory on $U_{k,N}(s)$, see, for example, \cite{ström} and Chapter 13 of \cite{MR3675870}.
\begin{proposition}\label{mellin}
The Mellin transform of $\varphi \in U_{0,N}(s)$ with Fourier expansion as in Proposition \ref{bo} at the cusp $\infty$ is given by
\begin{align*}
&\Omega(r) = \int\limits_{0}^{\infty}\varphi(iy)y^r\frac{dy}{y}\\
&= \pi^{-r-1/2}\Gamma\bigg(\frac{(r+1/2)+(s-1/2)}{2}\bigg)\Gamma\bigg(\frac{(r+1/2)-(s-1/2)}{2}\bigg)L\big(\varphi, r\big),
\end{align*}
where $L\big(\varphi, r\big)$ is the L-function of $\varphi$ which is defined by
\begin{align*}
L(\varphi, r) =\sum\limits_{n\geq 1}\frac{a(n)}{n^r}.
\end{align*}
It converges for $r > 1$. The function $\Omega(r)$ has an analytic continuation to all $r \in \mathbb{C}$, which we denote by $\Omega^*$ and which fulfills $\Omega^*(1-r) = \Omega^*(r)$.
\end{proposition}

\subsection{Vector Valued Maass Cusp Forms}
Let $\tau = u + iv \in \mathbb{H}$. Let $(V, q)$ be a rational quadratic space of signature $(b^+, b^-)$ and $ L \subset V$ an even lattice. We denote the group ring of $L'/L$ by $\mathbb{C}[L'/L]$ with basis vectors denoted by $\mathfrak{e}_h$ for $h\in L'/L$ and equip it with a scalar product $\langle \mathfrak{e}_h, \mathfrak{e}_{h'} \rangle = \delta_{h,h'}$, which is anti-linear in the second component.\\
We denote it by $\tilde{\Gamma} := \textup{Mp}_2(\mathbb{Z})$ the integral metaplectic group consisting of pairs $(M, \phi)$, with $M = \begin{psmallmatrix}a&b\\c&d \end{psmallmatrix} \in \textup{SL}_2(\mathbb{Z})$ and $\phi : \mathbb{H} \rightarrow \mathbb{C}$ holomorphic and such that $\phi(\tau)^2 = c\tau +d$. The Weil representation $\rho_L$ of $\tilde{\Gamma}$ is a unitary representation defined via its action on the generators $T =  \left(\begin{psmallmatrix}1&1\\0&1 \end{psmallmatrix}, 1\right)$ and $S = \left(\begin{psmallmatrix}0&-1\\1&0 \end{psmallmatrix}, \sqrt{z}\right)$ of $\tilde{\Gamma}$
\begin{align*}
&\rho_L(T)\mathfrak{e}_h = e^{2\pi i q(h)}\mathfrak{e}_h \\
&\rho_L(S)\mathfrak{e}_h = \frac{\sqrt{i}^{b^- - b^+}}{\sqrt{|L'/L|}}\sum\limits_{h' \in L'/L}e^{-2\pi i (h, h')}\mathfrak{e}_{h'}.
\end{align*}
The dual Weil representation we denote by $\rho_L^*$. We have $\rho_L^* = \overline{\rho}_L = \rho_{L^-}$.\\
Let $k \in \frac{1}{2}\mathbb{Z}$. Then, for any $(M, \phi) \in \textup{Mp}_2(\mathbb{Z})$ and any function $f: \mathbb{H} \rightarrow \mathbb{C}[L'/L]$, we define the weight $k$ \textbf{Petersson slash operator} by
\begin{align*}
(f|_{k, \rho_L}(M, \phi))(\tau) := \phi(z)^{-2k}\rho_L(M, \phi)^{-1}f(M\tau).
\end{align*}
We set $\tilde{\Gamma} = \textup{Mp}_2(\mathbb{Z})$. 
\begin{definition}\label{kakadu2}
A \textbf{Maass cusp form} of weight $k$ for $\rho_L$ is a real analytic function $F: \mathbb{H} \rightarrow \mathbb{C}[L'/L]$ which satisfies the following conditions:
\begin{itemize}
\item[(i)]$\Delta_kF = \lambda F$ for some $\lambda = (s - k/2) (1-s-k/2) \in \mathbb{C}$,
\item[(ii)] $(F|_{k, \rho_L}(M, \phi))(\tau) = F(\tau)$ for all $(M, \phi) \in \tilde{\Gamma}$,
\item[(iii)] There exists $C > 0$ such that $F(\tau) = \mathcal{O}(v^C)$ for $v \rightarrow \infty$,%muss das in jeder Spitze?,
\item[(iv)] The expression
\begin{align*}
\int\limits_0^1 F(u+iv)du,
\end{align*}
is equal to $0$ for all $v > 0$.
\end{itemize}
\end{definition}
We denote the finite dimensional space of all  Maass cusp forms of weight $k$ for $\rho_L$ and $\tilde{\Gamma}$ with fixed spectral point $s$ by $U_{k,\rho_L}(s)$.
\begin{proposition}
The Fourier series of a Maass cusp form $F(\tau) \in U_{k,\rho_L}(s)$ is given by
\begin{align}
\sum\limits_{h \in L'/L}\bigg[\sum\limits_{\substack {n \in \mathbb{Z} + q(h) \\n \neq 0}}c(h, n)(4\pi |n|v)^{-\frac{k}{2}}W_{\textup{sign}(n)\frac{k}{2}, s- \frac{1}{2}}(4\pi |n|v)e(nu)\bigg]\mathfrak{e}_h.
\end{align}
\end{proposition}
There is also a Hecke theory for vector valued Maass forms, see \cite{MR2574974} and \cite{str}. We will skip the technical details.
\begin{definition}
For the specific lattice $L$ as given in Section \ref{sig21}, $F(\tau) \in U_{k, \rho}(s)$ with $\rho = \rho_L$ or $\rho = \overline{\rho}_L$, we define the $p$-th Hecke operator for a prime $p$ with $(p,N) = 1$ by
\begin{align*}
    &F|_{k,L}T_{p^2}\\
    & = \sum\limits_{h \in L'/L}\big[\sum\limits_{\substack {n \in \mathbb{Z} + q(h) \\n \neq 0}}b(h, n)(4\pi |n|v)^{-\frac{k}{2}}W_{\textup{sign}(n)\frac{k}{2}, s- \frac{1}{2}}(4\pi |n|v)e(nu)\big]\mathfrak{e}_h,
\end{align*}
with
\begin{align*}
    b(h,n) = c(ph, p^2n) + p^{k-3/2}\begin{pmatrix}\frac{4\sigma Nn}{p}\end{pmatrix}c(h,n) + p^{2k-2}c(h/p, n/p^2),
\end{align*}
where $\sigma = 1$ if $\rho = \rho_L$ and $\sigma = -1$ if $\rho  = \overline{\rho}_L$.
\end{definition}
\begin{proposition}
    For $(m, N)=1$, the Hecke operators $T_{m^2}$ (the $T_{p^2}$ with $(p,N) =1))$ generate a commutative subalgebra of $\operatorname{End}(U_{k,\rho}(s))$. They are self-adjoint with respect to the Petersson inner product.
\end{proposition}
\begin{proof}
    This follows as in \cite{MR2574974}, Th. 4.12.
\end{proof}
\begin{corollary}\label{basis}
    The space $U_{k,\rho}(s)$ has a basis of simultaneous eigenforms for all $T_{n^2}$ with $(n, N) =1$.
\end{corollary}

\subsection{Differential Operators}\label{diffop}
We define the \textbf{Maass lowering operator} by $L_k := -2iy^2\frac{\partial}{\partial \overline{z}}$ and the \textbf{Maass raising operator} by $R_k :=  2i \frac{\partial}{\partial z}+ ky^{-1}$.
For a smooth function $f: \mathbb{H} \rightarrow \mathbb{C}[L'/L]$ and $(M, \phi) \in \textup{Mp}_2(\mathbb{Z})$, one has the following:
\begin{align*}
&R_k(f|_{k,\rho_L}(M,\phi)) = (R_kf)|_{k+2, \rho_L}(M,\phi),\\
&L_k(f|_{k,\rho_L}(M,\phi)) = (L_kf)|_{k-2, \rho_L}(M,\phi),\\
\end{align*}
These relations also hold for scalar valued functions and the slash operator for $\textup{SL}_2(\mathbb{Z})$.
For the hyperbolic Laplacian we have
\begin{displaymath}
\Delta_k = -\xi_{2-k}\xi_k = -L_{k+2}R_k - k = - R_{k-2}L_k.
\end{displaymath}
In particular, for $\Delta_k\varphi = \lambda\varphi = (s-k/2)(1-s-k/2)\varphi$, we have
\begin{align*}
&\Delta_{k-2}L_k\varphi = (\lambda - k+2)L_k\varphi = \left(s-\frac{k-2}{2}\right) \left(1-s-\frac{k-2}{2}\right)L_k\varphi,\\
&\Delta_{k+2}R_k\varphi = (\lambda +k)R_k\varphi = \left(s-\frac{k+2}{2}\right) \left(1-s-\frac{k+2}{2}\right)R_k\varphi.
\end{align*}

\begin{lemma}\label{der}
For any $n \in \mathbb{Q}$, we set
\begin{align*}
f_n(z) = a(n)(4\pi |n|y)^{-\frac{k}{2}}W_{\textup{sign}(n)\frac{k}{2}, s-\frac{1}{2}}(4\pi |n|y)e(nx).
\end{align*}
Then for the lowering operator we have
\begin{align*}
&L_kf_n(z) \\
&= \begin{cases}
 -a(n)(4\pi |n|)^{-\frac{k}{2}}y^{-\frac{(k-2)}{2}}W_{\textup{sign}(n)\frac{(k-2)}{2},s-\frac{1}{2}}(4\pi|n|y) e(nx), & \textup{ for } n <0,\\
 (s-k/2)(1-s-k/2)a(n)(4\pi |n|)^{-\frac{k}{2}}y^{-\frac{(k-2)}{2}}\\
\times W_{\textup{sign}(n)\frac{(k-2)}{2},s-\frac{1}{2}}(4\pi|n|y) e(nx), & \textup{ for } n > 0.
\end{cases}
\end{align*}
And for the raising operator we have
\begin{align*}
&R_kf_n(z) \\
&= \begin{cases}
 (s+k/2)(1-s+k/2)a(n)(4\pi |n|)^{-\frac{k}{2}}y^{-\frac{(k+2)}{2}}\\
\times W_{\textup{sign}(n)\frac{(k+2)}{2},s-\frac{1}{2}}(4\pi|n|y) e(nx), & \textup{ for } n < 0, \\
 -a(n)(4\pi |n|)^{-\frac{k}{2}}y^{-\frac{(k+2)}{2}}W_{\textup{sign}(n)\frac{(k+2)}{2},s-\frac{1}{2}}(4\pi|n|y) e(nx), & \textup{ for } n > 0.
\end{cases}
\end{align*}
\end{lemma}
\begin{proof}
This follows  by a standard computation involving the chain and product rule and the identities  on page 302 of \cite{MR232968}.
\end{proof}

\begin{proposition}\label{mc}
Let $\Gamma$ be a congruence subgroup \textup{(}$L$ an even lattice\textup{)}. The differential operators $L_k$ and $R_k$ map Maass cusp forms for $\Gamma$ \textup{(}$\tilde{\Gamma}$ and $\rho_L$\textup{)} to Maass cusp forms for $\Gamma$ \textup{(}$\tilde{\Gamma}$ and $\rho_L$\textup{)}.
\end{proposition}
\begin{proof}
This follows by a straight forward computation using the information above and Lemma \ref{der}
\end{proof}
We want to further characterize scalar valued Maass cusp forms. Therefore we consider their behaviour under a certain involution.
\begin{definition}\label{T}
 Assume that $k \in 2\mathbb{Z}$ and $s \neq \ell/2$ for any $\ell \equiv k (\textup{mod }2)$. We define a linear operator $T$ by
\begin{align*}
T:  U_{k, N}(s) \rightarrow U_{-k,N}(s), \hspace{3mm} \varphi(z) \mapsto y^k\varphi(-\overline{z}).
\end{align*}
For $k > 0$, we apply to a Maass form $\varphi(z) \in U_{k,N}(s)$ the iterated lowering operator $L_{k}^{k} = L_{2-k}\circ L_{k-4}\circ... \circ L_k$ and normalize the result by the factor $\frac{\Gamma(s-\frac{k}{2})}{\Gamma(s+\frac{k}{2})}$. For $k < 0$, we do the same with the raising operator and normalize by $\frac{\Gamma(s+\frac{k}{2})}{\Gamma(s-\frac{k}{2})}$. We end up with a Maass cusp form in $U_{-k,N}(s)$, to which we now apply $T$. This defines an involution on the space $U_{k, N}(s)$. We call this involution $I(\varphi)$. For further details on this map, consider Section 4 of \cite{MR1923476}.
\end{definition}
\begin{definition}\label{ivl}
A Maass cusp form $\varphi(z) \in U_{k,N}(s)$ is called \textbf{even} if its eigenvalue under $I$ is $+1$ and \textbf{odd} if it is $-1$.
\end{definition}
\begin{remark}
From now on, when using the term ``even", we automatically imply that the map above is a proper involution. This happens automatically when $s \neq \frac{\ell}{2}$ for any $\ell \equiv k (\textup{mod }2)$. Note that this excludes holomorphic cusp forms other than $\varphi \equiv 0$. Furthermore, for any Maass cusp form with non-exceptional eigenvalue $\lambda \neq 0$, this condition is always fulfilled.
\end{remark}
Let $\varphi \in U_{k,N}(s)$ be an even Maass cusp form and denote its $n$-th Fourier coefficient by $a(n)$. Let $n > 0$, then we have
\begin{align}\label{evn}
a(-n) = \frac{\Gamma(s+\frac{k}{2})}{\Gamma(s-\frac{k}{2})} a(n).
\end{align}
Furthermore, it can be checked by a quick computation, that every Maass cusp form which fulfills this relation is even and that the image of an even Maass cusp form under $R_k$ and $L_k$ is even.\\
Hecke eigenfunctions map to Hecke eigenfunctions under $R_k$ and $L_k$, which can be checked via a direct computation on the Fourier expansions see \cite{str}. Furthermore the Hecke operators commute with the Laplacian. This holds true for the vector and the scalar valued case.

\subsection{Heegner Points and Cycle Integrals}
In this section, we take a look at special points and geodesics related $L'$. Therefore, we set for $h \in L'/L$ and $D \in \mathbb{Q}$
\begin{align*}
L_{D,h} := \{X \in L+h \mid q(X) = D\}.
\end{align*}
The group $\Gamma$ acts on $L_{D,h}$ with finitely many orbits if $D \neq 0$.
\begin{subsubsection}{Heegner points}\label{heeg}
 For $D <0$, we can split $L_{D,h}$ into the sets $L_{D,h}^+$ and $L_{D,h}^-$, where the first contains the elements $\begin{psmallmatrix}x_2&x_1\\x_3&-x_2 \end{psmallmatrix}\in L_{D,h}$ with $x_3 \geq 0$ and the latter the elements with $x_3 < 0$. The case $x_3 = 0$ does not occur for $D < 0$.\\ \\
Recall that $\mathbb{D}$ is the Grassmannian of negative definite lines in $V(\mathbb{R})$. For $X \in V(\mathbb{Q})$ with $q(X) < 0$, we set
\begin{align*}
D_X := \textup{span}(X) \in \mathbb{D},
\end{align*}
and write $z_X$ for the corresponding point in $\mathbb{H}$. According to \cite{MR3353542}, Section 3, the point $z_X$ is given by $z_X = \frac{-b}{2a} + i\frac{\sqrt{|d|}}{2|a|}$ for $X = \begin{psmallmatrix}b&2c\\-2a&-b \end{psmallmatrix}$ and $q(X) = D = Nd < 0$. Furthermore, for $X \in V(\mathbb{Q})$ with $q(X) < 0$, the stabilizer $\Gamma_X$ of $X$ in $\Gamma$ is finite. Such points are called \textbf{CM-points} and the equivalence classes of such points under the action of $\Gamma$ are called \textbf{Heegner points}.
\end{subsubsection}
\begin{subsubsection}{Geodesics}
For $X \in V(\mathbb{Q})$ with $q(X) = D > 0$, we consider the set
\begin{align*}
C(X )= \{z \in \mathbb{D} \mid (X(z), X) = 0\} = \{z \in \mathbb{D} \mid z \perp X\}.
\end{align*}
The set $C(X)$ defines a geodesic in $\mathbb{D}$. See for example \cite{MR1244668} and \cite{MR1079646}.
Again after \cite{MR3353542}, for $X = \begin{psmallmatrix}b&2c\\-2a&-b \end{psmallmatrix}$, one obtains
\begin{align*}
C(X) := \{z \in \mathbb{D} \mid a|z|^2 + b\textup{Re}(z) + c = 0\}.
\end{align*}
The orientation of $C(X)$ is derived from that of the geodesic corresponding to $X_1 := \pm \begin{psmallmatrix}1&0\\0&-1\end{psmallmatrix}$, which is given by  the imaginary axis $\pm (0, i\infty)$ with indicated orientation. For any $X$ with $q(X) = D$, we find $g \in \textup{SL}_2(\mathbb{R})$ such that $g^{-1}X = \pm\sqrt{\frac{D}{N}} \begin{psmallmatrix}1&0\\0&-1\end{psmallmatrix}$. Since the action of $\textup{SL}_2(\mathbb{R})$ preserves the orientation, this automatically implies the orientation of $C(X)$.\\
We define a measure $dz_X$ in the following way: For $X_1 = \pm \sqrt{\frac{D}{N}} \begin{psmallmatrix}1&0\\0&-1\end{psmallmatrix}$, we define the measure $dz_{X_1} = \pm\frac{dz}{\sqrt{D}z}$ and set $dz_{g^{-1}X} = d(gz)_{X} $ for $g \in \textup{SL}_2(\mathbb{R})$ and $X \in V(\mathbb{Q})$. For $X =  \begin{psmallmatrix}b&a\\c&-b \end{psmallmatrix}$, this yields to
\begin{align*}
dz_X = \frac{2 \sqrt{N}dz}{Q_X(z)},
\end{align*}
with
\begin{align*}
Q_X(z) = N(-cz^2 + 2bz + a).
\end{align*}
Note that $Q_X(g z) = j(g, z)^{-2}Q_{g^{-1}X}(z)$ for all $g \in \textup{SL}_2(\mathbb{R})$.\\
The stabilizer $\overline{\Gamma}_X$ of $X$ in the image $\overline{\Gamma}$ of $\Gamma$ in $\textup{PSL}_2(\mathbb{Z})$ is either infinite cyclic or trivial. The first happens when $X^{\perp} \subset V$ is non-split over $\mathbb{Q}$, the second when $X^{\perp}$ is split. The last case arises if and only if $\frac{q(X)}{N} \in (\mathbb{Q}^{\times})^2$, see \cite{MR1930980}, Lemma 3.6 and Lemma 4.2.\\
We write $C_X$ for the image of $C(X)$ in $\Gamma\setminus\mathbb{D}$. If $\overline{\Gamma}_X$ is infinite cyclic, then $C_X$ is a closed geodesic in $\Gamma\setminus\mathbb{D}$. If it is trivial, then $C_X$ is an infinite geodesic. \\
For $X \in L'$ with $q(X) = D > 0$ and a Maass cusp form $f(z) \in U_{0,N}(s)$, we define the \textbf{cycle integral} corresponding to $X$ by
\begin{align*}
\int\limits_{C_X}f(z)dz_X.
\end{align*}
The integral converges due to the cuspidal properties of $f(z)$ and does not depend on the choice of $X \in \Gamma\setminus L_{D,h}$, since for any $g \in \textup{SL}_2(\mathbb{R})$, we have
\begin{align*}
\int\limits_{C_X}f(z)dz_X = \int\limits_{C_{g^{-1}X}}f(gz)dz_{g^{-1}X}.
\end{align*}
For a a Maass cusp form $f$ of weight $2m, m \in \mathbb{Z}$ and $C_Q$ a geodesic corresponding to a quadratic form $Q$ of discriminant $D > 0$, define
\begin{align*}
Cy(f, Q) = D^{\frac{1-m}{2}}\int\limits_{C_Q}f(z)Q(z,1)^{m-1}dz.
\end{align*}
This definition generalizes the definition of the cycle integral as above to higher weights by considering $Cy(f,X) = Cy(f, Q_X)$.
\end{subsubsection}

\subsection{Twisted Traces}
\begin{definition}
Let $\mathfrak{r} \in \mathbb{Z}$ and $\Delta \equiv \mathfrak{r}^2 (4N)$ be a fundamental discriminant and $D \in \mathbb{Q}$, $D < 0$. Let $h \in L'/L$. We define the \textbf{twisted trace} (of CM-values) of index $D$ and $h$ of a Maass cusp form $\varphi(z) \in U_{0,N}(s)$ by
\begin{align*}
Tr_{\Delta, \mathfrak{r}}(\varphi, h, D) = \sum\limits_{X\in \Gamma\setminus L_{|\Delta|D,\mathfrak{r}h}}\chi_{\Delta}(X)\frac{1}{|\overline{\Gamma}_X|}\varphi(z_X).
\end{align*}
We additionally define the terms
\begin{align*}
Tr^+_{\Delta, \mathfrak{r}}(\varphi, h, D) = \sum\limits_{X\in \Gamma\setminus L_{|\Delta|D,\mathfrak{r}h}^+}\chi_{\Delta}(X)\frac{1}{|\overline{\Gamma}_X|}\varphi(z_X),
\end{align*}
and $Tr^-_{\Delta, \mathfrak{r}}(\varphi, h, D)$ accordingly. Furthermore, we define
\begin{align*}
Tr^{\pm}_{\Delta, \mathfrak{r}}(\varphi, h, D) = Tr^+_{\Delta, \mathfrak{r}}(\varphi, h, D) - Tr^-_{\Delta, \mathfrak{r}}(\varphi, h, D).
\end{align*}
\end{definition}
\begin{definition}
For $\Delta \equiv \mathfrak{r}^2 (4N)$ a fundamental discriminant, $0 < D \in \mathbb{Q}$ and $h \in L'/L$, we define for a Maass cusp form $\varphi \in U_{0,N}(s)$ the expression
\begin{align*}
Tr_{\Delta, \mathfrak{r}}(\varphi, h, D) =  \sum\limits_{X\in \Gamma\setminus L_{|\Delta|D,\mathfrak{r}h}}\chi_{\Delta}(X)\int\limits_{C_X}\varphi(z)dz_X.
\end{align*}
We generalize this to all $k \in \mathbb{Z}$, $\varphi \in U_{k,N}(s)$ by defining the \textbf{twisted trace} (of geodesics) of weight $k$ for $D > 0$ by
\begin{align*}
Tr_{\Delta, \mathfrak{r},k}(\varphi, h, D) = \sum\limits_{X\in \Gamma\setminus L_{|\Delta|D,\mathfrak{r}h}}\chi_{\Delta}(X)Cy(\varphi,X).
\end{align*}
Note that for $k = 0$, both definitions agree up to a factor $\sqrt{|\Delta|D}$.
\end{definition}

\section{The Siegel and Millson Theta Functions}
Consider the quadratic space $(V,q)$ and the lattice $L$ as defined in Section \ref{sig21}. Recall that the congruence subgroup $\Gamma = \Gamma_0(N)$ acts isometrically on $V$, maps $L$ to itself and acts trivially on $L'/L$. Now, for $z = x + iy \in \mathbb{H}$ and $X = \begin{psmallmatrix}x_2 & x_1\\ x_3 & -x_2 \end{psmallmatrix} \in V(\mathbb{R})$, we define
\begin{align*}
&p_z(X) := \sqrt{2}(X, X_Z(z)) = \frac{\sqrt{N}}{y}(x_3|z|^2 -2x_2x -x_1), \\
&R(X, z) := (X,X) + \frac{1}{2}(X, X_Z(z))^2 = (X,X) + \frac{1}{2}p_z(X)^2.
\end{align*}
\begin{lemma}\label{identities}
For $\gamma \in \textup{SL}_2(\mathbb{Z})$, we have
\begin{align*}
p_{\gamma z}(X) = p_z(\gamma^{-1}X) \hspace{2mm} \textup{  and } \hspace{2mm} R(X, \gamma z) = R(\gamma^{-1}X, z).
\end{align*}
\end{lemma}
The quantity $R(X, z)$ is non-negative and vanishes exactly when $X \in \mathbb{R}X_Z(z)$. 
\begin{definition}\label{O}
For $\tau = u + iv \in \mathbb{H}$ and $z$ as above, we define the Schwartz functions
\begin{align*}
\varphi(z, \tau, X) := \sqrt{v}e^{-2\pi v R(X,z)}e^{2\pi i q(X) \overline{\tau}},
\end{align*}
and
\begin{align*}
\psi(z, \tau, X) := v^{3/2}p_z(X)e^{-2\pi v R(X,z)}e^{2\pi i q(X) \overline{\tau}}.
\end{align*}
\end{definition}

\begin{lemma}\label{lapla}
Let $\Delta_{1/2, \tau}$ be the weight $1/2$ Laplace operator in $\tau$ and $\Delta_{0, z}$ the weight $0$ Laplace operator in $z$, then we have the following identities:
\begin{align*}
\Delta_{1/2, \tau}\varphi(z, \tau, X) = \frac{1}{4}\Delta_{0,z}\varphi(z, \tau, X)
\end{align*}
and
\begin{align*}
\Delta_{1/2, \tau}\psi(z, \tau, X) = \frac{1}{4}\Delta_{0,z}\psi(z, \tau, X) - \frac{1}{2}\psi(z, \tau, X).
\end{align*}
\end{lemma}
\begin{proof}
For the first equation, see, for example, \cite{MR0389772}, Equation (2.10) of \cite{MR1244668}  or Section 5 of \cite{MR3353542}. For the second equation, consider \cite{craw}, Proposition 4.3.3.
\end{proof}
\begin{definition}The \textbf{Siegel theta function} for the lattice $L$ is defined by
\begin{align*}
\Theta_{S, L}(z, \tau) := \sum\limits_{h \in L'/L}\sum\limits_{X \in h+L}\varphi(z, \tau, X)\mathfrak{e}_h,
\end{align*}
and the \textbf{Millson theta function} for the lattice $L$ is defined by
\begin{align*}
\Theta_{M, L}(z, \tau) := \sum\limits_{h \in L'/L}\sum\limits_{X \in h+L}\psi(z, \tau, X)\mathfrak{e}_h.
\end{align*}
\end{definition}
\begin{proposition}
For the theta functions, we have the following: 
\begin{itemize}
\item[1.)] The Siegel theta function is an automorphic form of weight $0$ for $\Gamma_0(N)$ in $z$ and an automorphic form of weight $1/2$ for $\tilde{\Gamma}$ and $\rho_L$ in $\tau$.
\item[2.)] The Millson theta function is an automorphic form of weight $0$ for $\Gamma_0(N)$ in $z$ and an automorphic form of weight $-1/2$ for $\tilde{\Gamma}$ and $\rho_L$ in $\tau$.
\end{itemize}
\end{proposition}
\begin{proof}
The transformation behaviour in $z$ can be deduced from the identities in Lemma \ref{identities} and the fact that $\Gamma_0(N)$ acts trivially on $L'/L$.\\
For the transformation behaviour in $\tau$, we refer to Theorem 4.1 from \cite{MR1625724} and use the fact that $\textup{Im}(M\tau) = \frac{\textup{Im}(\tau)}{|c\tau +d|^2}$.\\
It is easy to see that $\Theta_{S, L}(z, \tau)$ is the Siegel theta function attached to the lattice $L$ (modified by a power of $\textup{Im}(\tau)$). For $\Theta_{M, L}(z, \tau)$, consider the isometry $\sigma$ from Lemma \ref{isomm}. Under this isometry, $p_z(X)$ is mapped to $\sqrt{2}x_3$ and thus, to a homogeneous polynomial of degree $(0,1)$. Furthermore, the polynomial $p_z(X)$ is harmonic. The rest now follows from Theorem 4.1 of \cite{MR1625724}.
\end{proof}

\subsection{Partial Poisson Summation}\label{43}
The general method we apply follows \cite{MR1625724}. For further reference, see \cite{MR2726107}, \cite{craw} and \cite{hoe}.
Throughout this section, we identify $z \in \mathbb{H}$ with $ z := \mathbb{R}X_2(z) + \mathbb{R}X_1(z) \in \mathbb{D}^{\perp}$.
\begin{definition}
Let $\ell \in L$ be a primitive isotropic vector and $\ell' \in L'$ such that $(\ell, \ell') = 1$. We define the lattice $K \subset L$ as
\begin{align*}
K := L \cap \ell'^{\perp} \cap \ell^{\perp}.
\end{align*}
The lattice $K$ is one-dimensional and positive definite. The space $K\otimes_{\mathbb{Z}}\mathbb{R}$ has signature $(1,0)$.
\end{definition}
For the sake of simplicity, we assume that $(\ell, L) = \mathbb{Z}$, so that we can choose the isotropic vector $\ell'$ from $L$. The lattice $L$ then splits into
\begin{align*}
L = K \oplus \mathbb{Z}\ell' \oplus \mathbb{Z}\ell,
\end{align*}
and we have $K'/K \cong L'/L$ and
\begin{align*}
V(\mathbb{R}) = L \otimes_{\mathbb{Z}}\mathbb{R} = (K\otimes_{\mathbb{Z}}\mathbb{R}) \oplus \mathbb{R}\ell' \oplus \mathbb{R}\ell .
\end{align*}
\begin{remark}
If the level $N$ of $L$ is square-free, every primitive isotropic vector $\ell \in L$ fulfills $(\ell, L) = \mathbb{Z}$. Throughout our computations, we restrict ourselves to square-free levels. %For general $N$ some modifications in the formulas might be necessary, but the general results of this section still hold.
\end{remark}
For the dual lattices $K'$ and $L'$, we have
\begin{align*}
K' = L' \cap \ell'^{\perp} \cap\ell^{\perp},
\hspace{2mm} \textup{ and }\hspace{2mm} L' = K' \oplus \mathbb{Z}\ell' \oplus \mathbb{Z}\ell.
\end{align*}
Let $\lambda_K$ be the projection of $\lambda \in V(\mathbb{R})$ onto $K\otimes \mathbb{R}$. For any $\lambda \in V(\mathbb{R})$, denote by $\lambda_z$ and $\lambda_{z^{\perp}}$ the orthogonal projections of $\lambda$ onto $z \in \mathbb{D}^{\perp}$  and $z^{\perp}$. Furthermore, denote by $w$ the orthogonal complement of $\ell_{z}$ in $z$. Then we have $z = \mathbb{R}\ell_{z}\oplus w$ and
\begin{align*}
V(\mathbb{R}) = z \oplus z^{\perp} = w \oplus \mathbb{R}\ell_z \oplus \mathbb{R}\ell_{z^{\perp}}.
\end{align*}
We define the vector
\begin{align*}
\mu = \mu(z) = \mu(z^{\perp}) := -\ell' + \frac{\ell_z}{2\ell_z^2}+\frac{\ell_{z^{\perp}}}{2\ell_{z^{\perp}}^2} \in V(\mathbb{R})\cap \ell^{\perp}.
\end{align*}
Note that the Grassmannian $\mathbb{D}_{K\otimes\mathbb{R}}$ only consists of a single point, which is why we omit the variable $z$ in the following construction.
\begin{remark}\label{papagai}
For $h = k + m\ell' + n\ell$ and $h' = k' + m'\ell' + n'\ell$ two vectors in $L' = K' \oplus \mathbb{Z}\ell' \oplus \mathbb{Z}\ell$, we have that $h = h'$ in $L'/L$ if and only if $k = k'$ in $K'/K$. Thus, for any $h \in L'/L$, we can always pick $h = k + 0\ell' +0\ell \in K'/K$. This means any sum over $L'/L$ can be interpreted as a sum over $K'/K$.
\end{remark}
\begin{definition}
Let $a, b \in \mathbb{Z}$, $\mu \in K\otimes_{\mathbb{Z}}\mathbb{R}$ and $h \in K'/K$. Then for $\tau \in \mathbb{H}$, we define a scalar valued twisted theta function
\begin{align*}
\xi_{h, \Delta}(\tau, \mu, a, b) := \sum\limits_{\lambda \in K + \mathfrak{r}h}\sum\limits_{\substack{j(\Delta)\\ q(\lambda - b\ell'+j\ell) \equiv \Delta q(h)(\Delta)}}\chi_{\Delta}(\lambda- b\ell'+j\ell)e(-aj/|\Delta|)\\
\times e\bigg(\frac{(\lambda + b\mu)^2}{2}\frac{\tau}{|\Delta|} - \frac{1}{|\Delta|}\bigg(\lambda + \frac{b\mu}{2}, a\mu\bigg)\bigg),
\end{align*}
and the corresponding $\mathbb{C}[K'/K]$-valued theta function
\begin{align*}
\Xi_{\Delta}(\tau, \mu, a, b) := \sum\limits_{h \in K'/K}\xi_h(\tau, \mu, a, b)\mathfrak{e}_h.
\end{align*}
\end{definition}

From now on, unless necessary, we will omit the index $\Delta$ from the theta function $\Xi_{\Delta}$ for the sake of simplicity. In the following, we will list a few properties of $\Xi(\tau, \mu, a, b)$ without proof.\\
\begin{proposition}[Theorem 4.4 of \cite{MR2726107}]  \label{poisson}
For $(M,\phi)\in \textup{Mp}_2(\mathbb{Z})$, where $M = \begin{psmallmatrix}\alpha & \beta\\ \gamma & \delta \end{psmallmatrix}$, we have
\begin{align*}
\Xi(M\tau, \mu, \alpha a + \beta b, \gamma a + \delta b)= \phi(\tau)\tilde{\rho}_K(M, \phi)\Xi(\tau, \mu, a, b),
\end{align*}
where $\tilde{\rho}_K$ is $\rho_K$ in case $\Delta > 0$ and $\overline{\rho}_K$ if $\Delta < 0$.
\end{proposition}
The following identities are crucial for computing the Fourier coefficients of the Shimura lift. See Lemma 4.7 and Theorem 4.8 of \cite{MR2726107} for the assertions concerning the Siegel theta function and Theorem 5.4.3 and 5.4.5 of \cite{craw} for the Millson theta function. Note that, in our case, the polynomial in the Millson theta function is $\sqrt{2/|\Delta|}$ times the one Crawford defines.
\begin{proposition}\label{thetah}
For the twisted Siegel theta function, we have
\begin{align*}
\theta_{S, \Delta, \mathfrak{r},h}(z, \tau) = \frac{1}{\sqrt{2|\Delta|\ell_z^2}}\sum\limits_{c,d \in \mathbb{Z}}e^{-\frac{\pi |c\tau +d|^2}{2|\Delta|v\ell_z^2}}\xi_h(\tau, \mu_K, d, -c),
\end{align*}
and for the Millson theta function, we have
\begin{align*}
\theta_{M,\Delta, \mathfrak{r},h}(z, \tau) = \frac{i}{2|\Delta|\ell_z^2}\sum\limits_{c,d \in \mathbb{Z}}(c\tau + d)e^{-\frac{\pi |c\tau +d|^2}{2|\Delta|v\ell_z^2}}\xi_h(\tau, \mu_K, d, -c).
\end{align*}
\end{proposition}
\begin{proposition} \label{platz}
For the Siegel theta function, we have that
\begin{align}
&\Theta_{S, \Delta}(z, \tau) = \frac{1}{\sqrt{2|\Delta|\ell_z^2}} \Xi(\tau, 0, 0, 0)\\
&+ \frac{1}{2\sqrt{2|\Delta|\ell_z^2}}\sum\limits_{n\geq 1}\sum\limits_{\gamma \in \tilde{\Gamma}_{\infty}\setminus \tilde{\Gamma}} \bigg[e^{\frac{-\pi n^2}{2|\Delta|\textup{Im}(\tau)\ell_z^2}}\Xi(\tau, \mu_K, n, 0)\bigg]\mid_{1/2, \tilde{\rho}_K}\gamma,\notag
\end{align}
and for the Millson theta function
\begin{align}
&\Theta_{M, \Delta}(z, \tau) = \frac{i}{4|\Delta|\ell_z^2} \sum\limits_{n\geq 1}\sum\limits_{\gamma \in \tilde{\Gamma}_{\infty}\setminus \tilde{\Gamma}}\bigg[e^{\frac{-\pi n^2}{2|\Delta|\textup{Im}(\tau)\ell_z^2}}\Xi(\tau, \mu_K, n, 0)\bigg]\mid_{-1/2, \tilde{\rho}_K}\gamma.
\end{align}
\end{proposition}
For our setup, since we mostly focus on the computations of the Fourier series in the cusp $\infty$, we take $\ell = \ell_{\infty} = \begin{psmallmatrix}0&1/N\\0&0 \end{psmallmatrix}$ and, consequently, we choose $\ell' = \begin{psmallmatrix}0&0\\1&0 \end{psmallmatrix}$. The lattice $K$ is then given by
\begin{align*}
K = \mathbb{Z}\begin{pmatrix}1&0\\0&-1 \end{pmatrix} \hspace{2mm} \textup{ and } \hspace{2mm} K' = \frac{1}{2N}\mathbb{Z}\begin{pmatrix}1&0\\0&-1 \end{pmatrix}.
\end{align*}
For $\lambda \in K\otimes\mathbb{R}$, we say $\lambda > 0$ if $\lambda = k \begin{psmallmatrix}1&0\\0&-1 \end{psmallmatrix}$ for $k > 0$. The lattice $K'$ can be identified with $\frac{1}{2N}\mathbb{Z}$ by identifying $\frac{m}{2N}\begin{psmallmatrix}1&0\\0&-1 \end{psmallmatrix}$ with $\frac{m}{2N}$, for $m \in \mathbb{Z}$. Thus, any sum over $K'$ can be read as a sum over $m \in \mathbb{Z}$ and we obtain for $\lambda \in K'$ that $q(\lambda) = \frac{m^2}{4N}$. For $z \in \mathbb{H}$ and under the above assumptions, we have, according to Section 5 of \cite{MR2726107}, that
\begin{align*}
&\ell_{z}^2 = \frac{1}{2N\textup{Im}(z)^2},\\
&(\lambda, \mu) = (\lambda, \textup{Re}(z)),\\
%t erklären\\
&\lambda^2/\ell_z^2 = (\lambda, \textup{Im}(z))^2.
\end{align*}
Furthermore, according to \cite{craw}, for $z = x + iy$, we have
\begin{align*}
&\mu = \begin{pmatrix}x& -x^2\\0&-x \end{pmatrix}  \hspace{2mm} \textup{ and } \hspace{2mm} \mu_K = \begin{pmatrix}x& 0\\0&-x \end{pmatrix},
\end{align*}
which gives us $(\lambda, \mu_K(i\eta)) = 0$ for $\eta >0$. Hence, we obtain the following lemma.
\begin{lemma}\label{ll}
We have
\begin{align*}
&\xi_h(\tau, \mu_K(i\eta), d, -c) \\
&=  \sum\limits_{\lambda \in K + \mathfrak{r}h}\sum\limits_{\substack{j(\Delta)\\ q(\lambda +c\ell'+j\ell) \equiv \Delta q(h)(\Delta)}}\chi_{\Delta}(\lambda +c\ell'+j\ell)e(-dj/|\Delta|)
 e\bigg(\frac{\lambda ^2}{2}\frac{\tau}{|\Delta|} \bigg),
\end{align*}
and
\begin{align*}
\xi_h(\tau, \mu_K(i\eta), n, 0) = \begin{pmatrix}\frac{\Delta}{n}\end{pmatrix}\overline{\varepsilon}|\Delta|^{1/2} \sum\limits_{\substack{\lambda \in K +\mathfrak{r}h\\ q(\lambda) = \Delta q(h) (\Delta)}} e\bigg(\frac{\lambda^2\tau}{2|\Delta|} \bigg).
\end{align*}
where $\varepsilon = 1$ if $\Delta > 0$ and $\varepsilon = i$ if $\Delta <0$. Furthermore, we set $\begin{pmatrix}\frac{\Delta}{0}\end{pmatrix}$ = 1 if $\Delta = 1$ and $\begin{pmatrix}\frac{\Delta}{0}\end{pmatrix} = 0$ otherwise.
\end{lemma}

\section{The Siegel and Millson Lift}\label{lifts}
Let $N \in \mathbb{N}_{>0}$  be square-free.
\begin{definition}
Let $k \in 4\mathbb{Z}$. Let $\varphi(z): \mathbb{H} \rightarrow \mathbb{C}$ be a function which is invariant under the weight $k$-slash operator for $\Gamma = \Gamma_0(N)$. We define (in case the beneath expression exists) the \textbf{twisted Siegel theta lift} of $\varphi(z)$ for the fundamental discriminant $\Delta$ by
\begin{align*}
I_{S, \Delta}(\tau, \varphi) = L_{1/2, \tau}^{-k/4}\int\limits_{\Gamma \setminus \mathbb{H}}R_{k, z}^{-k/2}\varphi(z)\Theta_{S, \Delta}(z,\tau)\frac{dxdy}{y^2}
\end{align*}
if $k < 0$ and by
\begin{align*}
I_{S, \Delta}(\tau, \varphi) = R_{1/2, \tau}^{k/4}\int\limits_{\Gamma \setminus \mathbb{H}}L_{k, z}^{k/2}\varphi(z)\Theta_{S, \Delta}(z,\tau)\frac{dxdy}{y^2}
\end{align*}
if $k \geq 0$.\\
Let now $k \in 4 \mathbb{Z} +2$ and $\varphi(z)$ defined as above. We define the \textbf{twisted Millson theta lift} of $\varphi(z)$ for $\Delta$ by
\begin{align*}
I_{M, \Delta}(\tau, \varphi) = L_{-1/2, \tau}^{(-k-2)/4}\int\limits_{\Gamma \setminus \mathbb{H}}R_{k, z}^{-k/2}\varphi(z)\Theta_{M, \Delta}(z,\tau)\frac{dxdy}{y^2}
\end{align*}
if $k < 0$ and by
\begin{align*}
I_{M, \Delta}(\tau, \varphi) = R_{-1/2, \tau}^{(k+2)/4}\int\limits_{\Gamma \setminus \mathbb{H}}L_{k, z}^{k/2}\varphi(z)\Theta_{M, \Delta}(z,\tau)\frac{dxdy}{y^2}
\end{align*}
if $k > 0$.
\end{definition}
\begin{proposition}\label{baum}
Let $s \neq \ell/2$ with $\ell \equiv k (\textup{mod } 2)$. For $\varphi \in U_{k,N}(s)$, the integrals $I_{S, \Delta}(\tau, \varphi)$ and $I_{M,\Delta}(\tau, \varphi)$ converge and lie in $U_{\frac{k+1}{2}, \tilde{\rho}_L}\left(\frac{s}{2}+\frac{1}{4}\right)$.
\end{proposition}
\begin{proof}
The integrals converge, since $\varphi(z)$ is of exponential decay towards the cusps, while, according to Proposition 5.4.6 of \cite{craw}, due to the action of the Atkin-Lehner involutions on automorphic forms of square-free level, at every cusp of $\Gamma$ and for $y \rightarrow \infty$, the Millson theta function lies in $\mathcal{O}\big(e^{-Cy^2}\big)$ and the Siegel theta function in $\mathcal{O}(y)$. The transformation behaviour in $\tau$ stems from that of the functions $\Theta_{S, \Delta}(z,\tau)$ and $\Theta_{M, \Delta}(z,\tau)$.\\
For the Laplace property, it follows for $k=0$ by Lemma \ref{lapla} and an argument as in \cite{MR1903920}, Lemma 4.3. 
The Laplace property in the general weight $k$ case follows with the fact that the raising and lowering operator map Laplace eigenfunctions to Laplace eigenfunctions and maintain the spectral point. The same arguments hold for the Millson theta lift.\\
It remains to show that the Siegel and Millson lift are again cusp forms. For the Millson theta function, we have by definition that $\psi(z, \tau, X) = \mathcal{O}(e^{-Cv}) \hspace{3mm} \textup{ for } v \rightarrow \infty$.
A standard argument involving dominated convergence then tells us that we can exchange limit and integration.\\
For the Siegel theta function, this follows as in \cite{MR1244668}. 
\end{proof}

\subsection{Fourier Expansions of Theta Lifts}\label{Siegelcoeff}
In this section we compute the Fourier expansion of the Siegel and the Millson theta lift, which we give in the following two theorems. 
\begin{theorem}\label{fünf}
Let $k \in 4\mathbb{Z}$ and $\varphi(z) \in U_{k,N}(s)$ be an even Maass cusp form. The Fourier expansion of the Siegel theta lift $I_{S, \Delta}(\tau, \varphi)$ in the cusp $\infty$ is given by
\begin{align*}
&I_{S, \Delta}(\tau, \varphi)\\
&\tiny{ = \sum\limits_{h \in L'/L}\sum\limits_{\substack {D \in \mathbb{Z}+\textup{sign}(\Delta)q(h)\\D \neq 0}}\tilde{C}(h, D)(4\pi |D|v)^{-\frac{k+1}{4}}W_{\textup{sign}(D)\frac{k+1}{4}, \frac{s}{2}- \frac{1}{4}}(4\pi |D|v)e(Du)\mathfrak{e}_h,}
\end{align*}
where for $k \geq 0$ we have
\begin{align*}
\tilde{C}(h, D)=\begin{cases}
(-1)^{\frac{k}{4}}(4\pi|D|)^{\frac{k}{4}}\frac{\sqrt{\pi}}{2}|D|^{-\frac{1}{2}}\frac{\Gamma\left(\frac{s}{2}+\frac{k+2}{4}\right)\Gamma\left(\frac{s}{2}\right)}{\Gamma\left(\frac{s}{2}+\frac{1}{2}\right)\Gamma\left(\frac{s}{2}-\frac{k}{4}\right)}Tr_{\Delta, \mathfrak{r}}\left(L_{k,z}^{k/2}\varphi, h, D\right), &D < 0\\ 
(-1)^{\frac{k}{4}}|D|^{-\frac{1}{2}}2^{\frac{k}{2}-1}(4\pi D)^{\frac{k}{4}}\frac{\Gamma\left(\frac{s}{2}+\frac{k+2}{4}\right)\Gamma\left(\frac{s}{2}\right)}{\Gamma\left(\frac{s}{2}+\frac{1}{2}\right)\Gamma\left(\frac{s}{2}-\frac{k}{4}\right)}\cdot Tr_{\Delta, \mathfrak{r},k}(\varphi, h, D), &D > 0.
\end{cases}
\end{align*}
and for $k < 0$ we have
\begin{align*}
\tilde{C}(h, D) =\begin{cases} 
(-1)^{\frac{k}{4}}(4\pi|D|)^{\frac{k}{4}}\frac{\sqrt{\pi}}{2}|D|^{-\frac{1}{2}}Tr_{\Delta, \mathfrak{r}}(R^{-k/2}\varphi, h, D),  & D < 0\\
(-1)^{\frac{k}{4}}(4\pi|D|)^{\frac{k}{4}}|D|^{-\frac{1}{2}}2^{\frac{k}{2}-1}\frac{\Gamma(s-k/2)}{\Gamma(s+k/2)}\cdot Tr_{\Delta, \mathfrak{r},k}(\varphi, h, D), & D > 0.
\end{cases}
\end{align*}
\end{theorem}

\begin{theorem}\label{sechs}
Let $k \in 4\mathbb{Z}+2$ and  $\varphi(z) \in U_{k,N}(s)$ be an even Maass cusp form. The Fourier expansion of the Millson theta lift $I_{M, \Delta}(\tau, \varphi)$ in the cusp $\infty$ is given by
\begin{align*}
&I_{M, \Delta}(\tau, \varphi)\\
&\tiny{= \sum\limits_{h \in L'/L}\sum\limits_{\substack {D \in \mathbb{Z} + \textup{sign}(\Delta)q(h) \\D \neq 0}}\tilde{C}(h, D)(4\pi |D|v)^{-\frac{k+1}{4}}W_{\textup{sign}(D)\frac{k+1}{4}, \frac{s}{2}- \frac{1}{4}}(4\pi |n|v)e(Du)\mathfrak{e}_h,}
\end{align*}
where for $k > 0$ we have
\begin{align*}
\tilde{C}(h, D) =
\begin{cases} (-1)^{\frac{k+2}{4}}\frac{(4\pi|D|)^{\frac{k+2}{4}}}{4\sqrt{\pi}|D|}\frac{\Gamma\left(\frac{s}{2}+\frac{k+2}{4}\right)\Gamma\left(\frac{s}{2}+\frac{1}{2}\right)}{\Gamma\left(\frac{s}{2}\right)\Gamma\left(\frac{s}{2}-\frac{k}{4}\right)} Tr^\pm_{\Delta, \mathfrak{r}}\left(L_{k,z}^{k/2}\varphi, h, D\right) &D < 0\\
 (-1)^{\frac{k+2}{4}}\frac{(4\pi|D|)^{\frac{k+2}{4}}2^{\frac{k+2}{2}}}{8\pi D}\frac{\Gamma\left(\frac{s}{2}+\frac{k+2}{4}\right)\Gamma\left(\frac{s}{2}+\frac{1}{2}\right)}{\Gamma\left(\frac{s}{2}\right)\Gamma\left(\frac{s}{2}-\frac{k}{4}\right)}\cdot Tr_{\Delta, \mathfrak{r},k}(\varphi, h, D), &D > 0.
\end{cases}
\end{align*}
and for $k < 0$ we have
\begin{align*}
\tilde{C}(h, D) = 
\begin{cases} (-1)^{-\frac{k+2}{4}}\frac{(4\pi|D|)^{\frac{k+2}{4}}}{4\sqrt{\pi}|D|}Tr^\pm_{\Delta, \mathfrak{r}}\left(R_{k,z}^{k/2}\varphi, h, D\right), &D < 0\\
(-1)^{-\frac{k+2}{4}}\frac{(4\pi|D|)^{\frac{k+2}{4}}}{8\pi D}2^{\frac{k+2}{2}}\frac{\Gamma(s-k/2)}{\Gamma(s+k/2)}\cdot Tr_{\Delta, \mathfrak{r},k}(\varphi, h, D), &D > 0.
\end{cases}
\end{align*}
\end{theorem}

\begin{proof}
We prove the results of Theorem \ref{fünf}, as Theorem \ref{sechs} is sufficiently similar. We first compute the integrals
\begin{align*}
I_{0,\Delta}(\tau, \varphi) = \int\limits_{\Gamma \setminus \mathbb{H}}R_{k, z}^{-k/2}\varphi(z)\Theta_{S, \Delta}(z,\tau)\frac{dxdy}{y^2}
\end{align*}
or
\begin{align*}
I_{0,\Delta}(\tau, \varphi) =\int\limits_{\Gamma \setminus \mathbb{H}}L_{k, z}^{k/2}\varphi(z)\Theta_{S, \Delta}(z,\tau)\frac{dxdy}{y^2}
\end{align*}
and apply the lowering or raising operator afterwards. We write down the computations only for $k \geq 0$, since the rest follows by the same methods. We first assume that $\Delta = 1$. The twisted result then follows as in \cite{MR3027941}.\\
We write $\vartheta_0(X, z, \tau)$ for $e^{-2\pi v R(X,z)}$ and define by $I(h,D)$ by the expression
\begin{align*}
I_{0,1}(\tau, \varphi) = \sum\limits_{h \in L'/L}\sum\limits_{D\in \mathbb{Q}}I(h,D)e(D\overline{\tau})\mathfrak{e}_h.
\end{align*}
We have
\begin{align}\label{frcf}
I(h,D) = \int\limits_{\Gamma \setminus \mathbb{H}}L_{k, z}^{k/2}\varphi(z)v^{1/2}\sum\limits_{X \in L_{D,h}}\vartheta_0(X, z, \tau)\frac{dxdy}{y^2}.
 \end{align}
The integral converges, due to the square-exponential decay of $R_{k, z}^{k/2}\varphi(z)$ and $L_{k, z}^{k/2}\varphi(z)$ towards the cusps of $\Gamma$ (see Proposition \ref{mc}). \\
 We end up with integral \eqref{frcf} given by
\begin{align}\label{bahnen}
&I(h,D) = \sum\limits_{X\in \Gamma\setminus L_{D,h}}\hspace{2mm}\int\limits_{\Gamma \setminus \mathbb{H}}L_{k, z}^{k/2}\varphi(z)v^{1/2}\sum\limits_{\gamma \in \Gamma_X\setminus \Gamma}\vartheta_0(\gamma X,z,\tau)\frac{dxdy}{y^2}\notag\\
& =  \sum\limits_{X\in \Gamma\setminus L_{D,h}}\hspace{2mm} \int\limits_{\Gamma_X \setminus \mathbb{H}}L_{k, z}^{k/2}\varphi(z)v^{1/2}\vartheta_0(X, z, \tau)\frac{dxdy}{y^2},
\end{align}
where the last identity is due to Lemma \ref{identities}.We consider the cases $D < 0$ and $D > 0$ separately.\\ \\
\textbf{Case 1}: $D < 0$.\\
For $D < 0$, the stabilizer group of every vector $X \in V$, with $q(X) = D$ in $\Gamma$ is finite and the expression in \eqref{bahnen} simplifies to
\begin{align*}
\sum\limits_{X \in \Gamma\setminus L_{D,h}}\frac{1}{|\overline{\Gamma}_X|} \int\limits_{ \mathbb{H}}L_{k, z}^{k/2}\varphi(z)v^{1/2}\vartheta_0(X, z, \tau)\frac{dxdy}{y^2}\\
=: \sum\limits_{X\in \Gamma\setminus L_{D,h}}\frac{1}{|\overline{\Gamma}_X|} I(X).
\end{align*}
%For every $z \in \mathbb{H}$, we find $g \in \textup{SL}_2(\mathbb{R})$ such that $gi = z$.
We use the approach of Katok and Sarnak in \cite{MR1244668} and write $I(X)$ as an integral over $\textup{SL}_2(\mathbb{R})$ in the following way
\begin{align*}
I(X) = \int\limits_{\textup{SL}_2(\mathbb{R})}L_{k, z}^{k/2}\varphi(gi)v^{1/2} \vartheta_0(X, gi, \tau) dg.
\end{align*}
Here, the Haar measure $dg$ is normalized so that the maximal compact open subgroup has volume 1. Note that we must only compute the integral $I(X)$ for $X \in L_{D, h}^+$, since for $X \in L_{D, h}^-$, we have that $-X \in L_{D, -h}^+$ and $\vartheta_0(-X, gi, \tau) = \vartheta_0(X, gi, \tau)$. Following Section 3 of \cite{MR1244668},
since for $D < 0$, the group $\textup{SL}_2(\mathbb{R})$ acts transitively on $\Big\{ \begin{psmallmatrix}x_2 & x_1\\ x_3 & -x_2 \end{psmallmatrix} \in \textup{M}_2(\mathbb{R}) \mid x_2^2 -x_1x_3 = D/N, x_3 > 0\Big\}$,  we find $g_0$ such that $g_0X_{D,0} := \begin{psmallmatrix}x_2 & x_1\\ x_3 & -x_2 \end{psmallmatrix} = X$ for all $x_1, x_2, x_3 \in \mathbb{Q}$. Here, we have set
\begin{align*}
X_{D,0} = -\frac{\sqrt{|D|}}{\sqrt{N}} \begin{pmatrix}0 & 1\\ -1 & 0 \end{pmatrix} = -\sqrt{2|D|}X_Z(i).
\end{align*}
We furthermore note that $g_0i$ is the Heegner point $z_X$ in $\mathbb{H}$, corresponding to $D_X$, see Section \ref{heeg} and Section \ref{sig21}.
For our integral, we use the invariance of the Haar measure and the properties of $R(X,z)$, as given in Lemma \ref{identities}, to obtain that
\begin{align*}
I(X) = \int\limits_{\textup{SL}_2(\mathbb{R})}L_{k, z}^{k/2}\varphi(g_0gi)v^{1/2}\vartheta_0(g^{-1}X_{D,0}, i, \tau) dg.
\end{align*}
By means of the Cartan decomposition, we have $SL_2(\mathbb{R}) = KA^+K$, where $K = \textup{SO}(2)$ and $A^+ = \bigg\{\alpha(a) = \begin{psmallmatrix}a&0\\0&a^{-1} \end{psmallmatrix}\mid a \geq 1\bigg\}$.\\ %Quelle
We set $\varphi_0(g) := \varphi(g_0gi)$. Using all of this, we obtain
\begin{align*}
&I(X) = 4\pi \int\limits_K\int\limits_{1}^{\infty}\int\limits_K L_{k, z}^{k/2} \varphi_0(k_1\alpha(a)k_2)v^{1/2}\\
&\times \vartheta_0((k_1\alpha(a)k_2)^{-1}X_{D,0}, i, \tau) \frac{a^2-a^{-2}}{2}dk_1\frac{da}{a}dk_2. %da fehlt 4\pi??
\end{align*}
Since the stabilizer of $X_{D,0}$ and $i$ is exactly given by $K= \textup{SO}(2)$, this is equal to
\begin{align*}
 4\pi v^{1/2}\int\limits_{1}^{\infty}\vartheta_0(\alpha(a)^{-1}X_{D,0}, i, \tau)\int\limits_K\int\limits_K L_{k, z}^{k/2}\varphi_0(k_1\alpha(a)^{-1}k_2)dk_1dk_2\frac{a^2-a^{-2}}{2}\frac{da}{a}.
\end{align*}
Now, the function $\varphi_0(k_1\alpha(a)k_2)$ is an eigenfunction of the Casimir operator with the same eigenvalue $\lambda$ as $\varphi(z)$ (see \cite{MR0803508}, Chapter X, §2). It is also $K$-biinvariant. %Casimir operator definieren!
Thus, the integral
\begin{align*}
f_1(g) = f_1(\alpha(a)) := \int\limits_K\int\limits_K L_{k, z}^{k/2}\varphi_0(k_1\alpha(a)k_2)dk_1dk_2,
\end{align*}
can be written as
\begin{align*}
f_1(g) = \varphi_0(1)\omega_{\lambda}(g) =  \varphi_0(1)\omega_{\lambda}(\alpha(a)) ,
\end{align*}
where $\omega_{\lambda}(g)$ is the unique spherical function of eigenvalue $\lambda$ as in Chapter IV, §3 of \cite{MR0803508}. 
By definition, we have $\varphi_0(1) =  \varphi(g_0i) = \varphi(z_X)$. Thus, we obtain
\begin{align*}
I(X) =  \bigg[ L_{k, z}^{k/2}\varphi(z_X)\bigg] Y_{\lambda}\bigg( \frac{\sqrt{|D|}}{\sqrt{N}}\bigg),
\end{align*}
with
\begin{align*}
Y_{\lambda}\bigg( \frac{\sqrt{|D|}}{\sqrt{N}}\bigg) = 4\pi v^{1/2}\int\limits_{1}^{\infty}\vartheta_0(\alpha(a)^{-1}X_{D,0}, i, \tau)\omega_{\lambda}(\alpha(a))\frac{a^2-a^{-2}}{2}\frac{da}{a}.
\end{align*}
The finction $\omega_{\lambda}(\alpha(a))$  only depends on the value $\frac{a^2+a^{-2}}{2}$, see Chapter V, §2 of \cite{MR0803508}. Thus, we write  $\omega_{\lambda}(\alpha(a))= \omega_{\lambda}\Big(\frac{a^2+a^{-2}}{2}\Big)$. We substitute $a = e^{r/2}$, then $x = \sinh(r)^2$. 
\begin{displaymath}
\vartheta_0(\sqrt{2}t\alpha(a)^{-1}X_Z(i), i, \tau) = e^{-4\pi v t^2 x}.
\end{displaymath}
The spherical function $\omega_{\lambda}$ is then a function in $\cosh(r) = \sqrt{1+x}$ and, as such, according to Chapter 9.2 of \cite{MR0499000}, is given by the Legendre function $P_{s-1}(\sqrt{1+x})$ of eigenvalue $\lambda = s(1-s)$.
 %evtl prüfen, ob das Polynom stimmt
Thus, we obtain 
\begin{align*}
Y_{\lambda}\bigg( \frac{\sqrt{|D|}}{\sqrt{N}}\bigg) = \pi v^{1/2}\int\limits_{0}^{\infty}e^{-4\pi v |D| x}(\sqrt{1+x})^{-1}P_{s-1}(\sqrt{1+x}) dx.
\end{align*}
From page 180, row 8 in \cite{MR0061695}, we see that the integral in the equation above is given by
\begin{align*}
(4\pi |D|v)^{-3/4} W_{-\frac{1}{4}, \frac{s}{2}-\frac{1}{4}}(4\pi |D|v)e^{2\pi |D|v}.
\end{align*}

\noindent \textbf{Case 2}: $D > 0$.\\
This time, we denote
\begin{align*}
J(X) := \int\limits_{\Gamma_X \setminus \mathbb{H}}L_{k, z}^{k/2}\varphi(z)v^{1/2}\vartheta_0(X, z, \tau)\frac{dxdy}{y^2}.
\end{align*}
Similar to before, we now find $h_0 \in \textup{SL}_2(\mathbb{R})$ such that $h_0X_{D,1} = \begin{psmallmatrix}x_2 & x_1 \\ x_3 & -x_2 \end{psmallmatrix}$ for all $x_1, x_2, x_3 \in \mathbb{Q}$, where
\begin{align*}
X_{D,1} = \frac{\sqrt{D}}{\sqrt{N}}\begin{pmatrix}1&0\\0&-1 \end{pmatrix}.
\end{align*}
We consider the integral over $h_0z$ instead of $z$ and use Lemma \ref{identities} to obtain
\begin{align*}
J(X)  %\int\limits_{\Gamma_X' \setminus \mathbb{H}}L_{k, z}^{k/2}\varphi(h_0z)v^{1/2} \vartheta_0(X, h_0z, \tau)\frac{dxdy}{y^2}\\
=  \int\limits_{\Gamma_X' \setminus \mathbb{H}}L_{k, z}^{k/2}\varphi(h_0z)v^{1/2}\vartheta_0(X_{D,1}, z, \tau)\frac{dxdy}{y^2},
\end{align*}
where $\Gamma_X' = h_0^{-1}\Gamma_Xh_0$. From this, we see that $\Gamma_X'$ is a discrete subgroup of the stabilizer of $X_{D,1}$ in $\textup{SL}_2(\mathbb{R})$ which is given by $\left\{\pm \begin{psmallmatrix}a & 0\\ 0 & a^{-1} \end{psmallmatrix} \biggm| a> 0 \right\}$.
There are two possibilities: Either $\Gamma_X'$ is infinite cyclic (if $D/N$ is not a square) or it is trivial (i.e.\ $\Gamma_X' = \{\pm 1\}$ if $D/N$ is  a square). We only consider the former case, since the latter is sufficiently similar. We obtain
\begin{align*}
\Gamma_X' =  \bigg\{\pm \begin{pmatrix}\varepsilon & 0\\ 0 & \varepsilon^{-1} \end{pmatrix}^n \biggm| n \in \mathbb{Z}\bigg\}
\end{align*}
for $\varepsilon > 1$. A fundamental domain of $\Gamma_X' \setminus \mathbb{H}$ is then given by $\{z \in \mathbb{H} \mid 1 \leq y < \varepsilon^2 \}$. We substitute $t = \frac{x}{y}$ in the integral over $x$ and obtain:
\begin{align*}
J(X) = v^{1/2}\int\limits_{-\infty}^{\infty}\biggl(\int\limits_{1}^{\varepsilon^2}L_{k, z}^{k/2}\varphi(h_0(y(t+i)))\frac{dy}{y}\biggr)e^{-4\pi v D (t^2+1)}dt.
\end{align*}
Following Katok and Sarnak again, %Katok Sarnak %POLYNOM MUSS RAUS!
we define
\begin{align*}
H(t) := \int\limits_{1}^{\varepsilon^2}L_{k, z}^{k/2}\varphi(h_0(y(t+i)))dk\frac{dy}{y}.
\end{align*}
We substitute $t=\frac{x}{y}$ into the Laplace differential equation and see that, since $L_{k, z}^{k/2}\varphi(h_0z)$ is an eigenfunction of the Laplace operator with eigenvalue $s(1-s)$, the function $H(t)$ solves
\begin{align}\label{eq}
(1+t^2)H''(t) + 2tH'(t) - s(1-s)H(t) = 0.
\end{align}
When substituting $z = -t^2$ and denoting $\tilde{H}(z) = H(t)$, this becomes a hypergeometric differential equation, namely
\begin{align*}
z(1-z)\tilde{H}''(z) + \left(\frac{1}{2} - \frac{3}{2}z\right)\tilde{H'}(z) - \frac{s(1-s)}{4}\tilde{H}(z) = 0,
\end{align*}
which, according to Chapter 2 of \cite{MR0029000}, around $z = 0$ has the solutions
\begin{align*}
\tilde{V}_{\lambda}(z) = F\Big(\frac{s}{2}, \frac{1}{2}-\frac{s}{2}, \frac{1}{2}; z \Big) \textup{ and } \tilde{U}_{\lambda}(z) = z^{1/2}F\Big(\frac{1}{2}+ \frac{s}{2}, 1-\frac{s}{2}, \frac{3}{2}; z \Big).
\end{align*}
Thus, we find that the unique even and odd solution of \eqref{eq} are given by
\begin{align*}
V_{\lambda}(t) = F\Big(\frac{s}{2}, \frac{1}{2}-\frac{s}{2}, \frac{1}{2}; -t^2 \Big) \textup{ and } U_{\lambda}(t) = tF\Big(\frac{1}{2}+ \frac{s}{2}, 1-\frac{s}{2}, \frac{3}{2}; -t^2 \Big),
\end{align*}
where $V(0) = U'(0) = 1$ and $U(0) = V'(0) = 0$. Subsequently, we can write
\begin{align}\label{dg}
H(t) = H(0)V_{\lambda}(g) + H'(0)U_{\lambda}(g).
\end{align}
The integral over the odd part of the input function vanishes, so we obtain
\begin{align*}
J(X) = H(0) v^{1/2}\int\limits_{-\infty}^{\infty}V_{\lambda}(t)e^{-4\pi vD(t^2 + 1)}dt.
\end{align*}
We see that
\begin{align*}
H(0) = \int\limits_{1}^{\varepsilon^2}L_{k,z}^{k/2}\varphi(h_0yi)\frac{dy}{y} = \sqrt{D}\int\limits_{c_X}L_{k,z}^{k/2}\varphi(z)dz_X, %Zykelintegralarbeit zitieren
\end{align*}
where the last equation is due to parametrising $C_X$ by $h_0iy$ for $y \in [1, \varepsilon^2]$, using $\frac{d(h_0yi)}{dy}= j(h_0, iy)^{-2}i$ and taking into account the identity
\begin{align*}
Q_X(h_0yi) = j(h_0, iy)^{-2} (2\sqrt{ND} iy).
\end{align*}
 It remains to understand the expression
\begin{align*}
 \int\limits_{-\infty}^{\infty}V_{\lambda}(t)e^{-4\pi vD(t^2 + 1)}dt
= e^{-4\pi vD} \int\limits_{-\infty}^{\infty}F\Big(\frac{s}{2}, \frac{1}{2}-\frac{s}{2}, \frac{1}{2}; -t^2 \Big)e^{-4\pi vDt^2}dt.
\end{align*}
We take into account the symmetry of the integral and substitute $t^2 = w$ to obtain that it is equal to
\begin{align*}
e^{-4\pi vD} \int\limits_{0}^{\infty}\frac{1}{\sqrt{w}}F\Big(\frac{s}{2}, \frac{1}{2}-\frac{s}{2}, \frac{1}{2}; -w \Big)e^{-4\pi vDw}dw,
\end{align*}
which, according to equation 3.37.1.2 of \cite{MR1162979}, is a Laplace transform given by
\begin{align*}
\sqrt{\pi}(4\pi D v)^{-3/4}W_{\frac{1}{4}, \frac{s}{2}-\frac{1}{4}}(4\pi Dv)e^{-2\pi Dv}.%Achtung Whittaker-Faktor stimmt nicht mit Buch uberein
\end{align*}
From this and case one, we see that the untwisted Siegel theta lift is given by
\begin{align*}
&I_{0,1}(\tau, \varphi) = \sum\limits_{h \in L'/L}\sum\limits_{D\in \mathbb{Q}}I(h,D)e(D\overline{\tau})\mathfrak{e}_h\\
&=  \sum\limits_{h \in L'/L}\big[\sum\limits_{\substack {D \in \mathbb{Q}\\D \neq 0}}C_0(h,D)(4\pi |D|v)^{-\frac{1}{4}}W_{\textup{sign}(D)\frac{1}{4}, \frac{s}{2}- \frac{1}{4}}(4\pi |D|v)e(Du)\big]\mathfrak{e}_h,
\end{align*}
with
\begin{align*}
C_0(h,D) =
\begin{cases}
\sqrt{\pi}\frac{|D|^{-1/2}}{2}Tr_{1,1}\left(L_{k,z}^{k/2}\varphi, h, D\right) & \textup{ if } D < 0,\\
2^{-1}Tr_{1,1}\left( L_{k,z}^{k/2}\varphi, h, D \right) & \textup{ if } D > 0.
\end{cases}
\end{align*}
Here, we have used that $Tr^-_{1,1}(f, h, D) = Tr^+_{1,1}(f, -h, D)$ for all $f\in U_{0,N}(s)$. Note as well that we already know that $I_S(\tau, \varphi)$ is a cusp form, due to Proposition \ref{baum}, so there is no $0$-th coefficient to consider. \\ \\
\noindent \textbf{Non-zero weights}:\\
Finally, we need to apply the raising (or lowering) operator to the  weight $1/2$ Fourier series
\begin{align}\label{for}
&I_{0, \Delta}(\tau, \varphi) \\
&\small{= \sum\limits_{h \in L'/L}\big[\sum\limits_{\substack {D \in \mathbb{Z}+\textup{sign}(\Delta)q(h) \\D \neq 0}}C(h, D)(4\pi |D|v)^{-\frac{1}{4}}W_{\textup{sign}(D)\frac{1}{4}, \frac{s}{2}- \frac{1}{4}}(4\pi |D|v)e(Du)\big]\mathfrak{e}_h.}
\end{align}
Denote the result of this procedure by
\begin{align*}
\tiny{\sum\limits_{h \in L'/L}\big[\sum\limits_{\substack {D \in \mathbb{Z}+\textup{sign}(\Delta)q(h) \\D \neq 0}}\tilde{C}(h, D)(4\pi |D|v)^{-\frac{k+1}{4}}W_{\textup{sign}(D)\frac{k+1}{4}, \frac{s}{2}- \frac{1}{4}}(4\pi |D|v)e(Du)\big]\mathfrak{e}_h,}
\end{align*}
then we apply Lemma \ref{der} and see that we have for $D > 0$ that
\begin{align*}
\tilde{C}(h, D) = (-1)^{k/4}(4\pi D)^{k/4}C(h,D),
\end{align*}
and for $D <0$ that
\begin{align*}
\tilde{C}(h, D) = (-1)^{k/2}(4\pi |D|)^{k/4}\frac{\Gamma\left(\frac{s}{2}+\frac{k+2}{4}\right)\Gamma\left(\frac{s}{2}\right)}{\Gamma\left(\frac{s}{2}+\frac{1}{2}\right)\Gamma\left(\frac{s}{2}-\frac{k}{4}\right)}C(h,D).
\end{align*}
We want to take a closer look at the positive Fourier coefficients. We only consider the untwisted ones, as the twisted analogue follows directly by \cite{MR3027941}. 
\begin{align*}
&Tr_{1,1,0}\left(L_{k,z}^{k/2}\varphi, h, D\right) = D^{1/2}Tr_{1,1}\left(L_{k,z}^{k/2}\varphi, h, D\right) = D^{1/2}\int\limits_{c(X)}L_k^{k/2}\varphi(z)dz_X\\
&= D^{1/2}\int\limits_{c(X)}L_{k,z}^{k/2}\varphi(z)\frac{2\sqrt{N}dz}{Q(z,1)} = Cy\left(L_{k,z}^{k/2}\varphi(z), Q\right),
\end{align*}
where $Q:= Q_X$ is the quadratic form corresponding to $X \in L'$. Note that $2\sqrt{N}D^{1/2}$ is exactly the square root of the discriminant corresponding to $Q$. We have
\begin{align*}
&Cy\left(L_k^{k/2}\varphi(z), Q\right) = Cy\left(L_4^2\left(L_k^{k/2-2}\varphi(z)\right), Q\right)\\
&= Cy\left(L_4^2\left(L_8^2\left( ... (L_k^{2}\varphi(z))\right)\right), Q\right).
\end{align*}
We now successively apply Theorem 1.1 of \cite{MR4122768}. Note that if the exponent of the lowering operator is $2$ in \cite{MR4122768}, the quantity $\ell$ from the same Theorem is equal to $-m$, where, in our case, we set $m_0 = -3$ and $m_j = m_{j-1}-2$, so that $m_1= - 5, ...,m_{\frac{k-4}{4}}= -\frac{k+2}{2}$. Note as well that $k$ is divisible by $4$. By Theorem 1.1 of \cite{MR4122768}, we have
\begin{align*}
&Cy\left(L_{2-2m_j}^2\left(L_{2-2(m_j-2)}^2\left(...  (L_k^{2}\varphi(z))\right)\right)\right) \\
&= -\lambda_{m_j}Cy\left(L_{2-2m_{j+1}}^2\left(...  (L_k^{2}\varphi(z))\right)\right),
\end{align*}
where $-\lambda_{m_j} = (s + m_j -1) (s - m_j)$ is $(-1)$-times the eigenvalue of\\ $L_{2-2m_{j+1}}^2(...  (L_k^{2}\varphi(z)))$.
We obtain that
\begin{align*}
Cy\left(L_k^{k/2}\varphi(z), Q\right)
%& = (s-2)(s+1)(s-4)(s+3)...(s-k/2)(s + k/2-1) \cdot Cy\left(\varphi(z), Q\right)\\
= 2^{k/2}\frac{\Gamma\left(\frac{s}{2}+\frac{k+2}{4}\right)\Gamma\left(\frac{s}{2}\right)}{\Gamma\left(\frac{s}{2}+\frac{1}{2}\right)\Gamma\left(\frac{s}{2}-\frac{k}{4}\right)}\cdot Cy\left(\varphi(z), Q\right).
\end{align*}
We have $Cy\left(\varphi(z), Q\right) = Tr_{1,1,k}\left(\varphi, h, D\right)$. Taking this into account, we arrive at the Fourier series given. \\
Now, for $k<0$, we need to exchange lowering and raising operator. This is straightforward for the negative Fourier coefficients. For the positive coefficients, by Lemma \ref{der}, the repeated application of the lowering operator to the Fourier series in \eqref{for} gives us that
\begin{align*}
&\tilde{C}(h, D) = (-1)^{k/4}(4\pi|D|)^{k/4} \frac{\Gamma\left(\frac{s}{2}-\frac{k}{4}\right)\Gamma\left(\frac{s}{2}+\frac{1}{2}\right)}{\Gamma\left(\frac{s}{2}\right)\Gamma\left(\frac{s}{2}+\frac{k+2}{4}\right)}C(h,D)\\
&= (-1)^{k/4}(4\pi|D|)^{k/4} 2^{k/2}\prod\limits_{j=0}^{-\frac{k+4}{4}}(s+2j)(s-2j-1)C(h,D).
\end{align*}
However, for the cycle integral $Cy\left(R^{-k/2}_{k,z}\varphi, Q \right)$, we have
\begin{align*}
Cy\left(R^{-k/2}_{k,z}\varphi, Q \right) = Cy\left(R^2_{-4,z}\left(R_{-8,z}^2\left(... (R^2_{k,z}\varphi)\right)\right), Q\right).
\end{align*}
We repeatedly apply Theorem 1.1 of \cite{MR4122768}, where, this time, we have $\ell = m_j-2$  with $m_0 = \frac{4+2}{2}$ and $m_j = m_{j-1} + 2$, so that $m_1= \frac{8+2}{2}, ..., m_{-\frac{k+4}{4}} = \frac{-k+2}{2}$. Applying Theorem 1.1 of \cite{MR4122768} in  for every $m_j$ gives us
\begin{align*}
&R_{2-2m_j, z}^2R_{2-2(m_j+2), z}^2( ...(R_{k,z}^2\varphi))\\
& = ((2m_j -2)-\lambda_{m_j}) R_{2-2(m_{j+1}), z}^2( ...(R_{k,z}^2\varphi)).
\end{align*}
Here, by $\lambda_{m_j}$ we denote the eigenvalue of $R_{2-2(m_{j+1}), z}^2( ...(R_{k,z}^2\varphi))$, a function of weight $2-2m_j$. However, $((2m_j -2)-\lambda_{m_j})$ is exactly the eigenvalue of $R_{2-2m_j, z}(R_{2-2(m_{j+1}), z}^2( ...(R_{k,z}^2\varphi)))$ times $(-1)$, which is given by $(s + 1-m_j)(s + m_j-2)$. We have $m_j = 3 + 2j$ for $ j \in \{0, ..., -\frac{k+4}{4}\}$, so we obtain
\begin{align*}
Cy\left(R^{-k/2}_{k,z}\varphi, Q \right) =  \prod\limits_{j=0}^{-\frac{k+4}{4}}(s + 2j +1)(s -2j -2)Cy(\varphi, Q).
\end{align*}
In summary, we arrive at
\begin{align*}
\tilde{C}(D, h) = (-1)^{k/4}|\Delta D|^{-1/2}(4\pi|D|)^{-k/4}2^{k/2}\frac{\Gamma(s-k/2)}{\Gamma(s+k/2)}Tr_{\Delta, \mathfrak{r},k}(\varphi, h, D).
\end{align*}
\end{proof}

\subsection{The Siegel and Millson Theta Lift for Switched Signature}
Due to the properties of the lowering and raising operator, we had to restrict ourselves to certain weights for each, i.e weight $k \in 4\mathbb{Z}$ for the Siegel theta lift and weight $k \in 4\mathbb{Z}+2$ for the Millson theta lift. We show how to get rid of this restriction and arrive at a formula covering all weights for both lifts, once the restricted result is known.
We consider the action of the operator $T$ from Definition \ref{T} on vector valued Maass cusp forms. It sends a Maass cusp form $\psi(\tau)$ to $v^k\psi(-\overline{\tau})$.

\begin{lemma}
For any $\psi(\tau) \in U_{k, \rho_L}(s)$, we have that $T\circ\psi := T\psi \in  U_{-k, \overline{\rho}_L}(s)$.
\end{lemma}
\begin{proof}
For $M =\begin{psmallmatrix}a&b \\ c & d \end{psmallmatrix}$, set $M' = \begin{psmallmatrix}a&-b \\ -c & d \end{psmallmatrix}$. Define an involution $v: \tilde{\Gamma} \rightarrow \tilde{\Gamma}$ with
\begin{align*}
v\left(M, \phi \right) = \left(M', \overline{\phi(-\overline{\tau})}\right).
\end{align*}
It can easily be checked that $v((M_1, \phi_1)(M_2, \phi_2)) = v((M_1, \phi))v((M_2, \phi))$.  We have $\rho_L( v(M,\phi)) = \overline{\rho}_L(M, \varphi)$, which can be verified by a direct computation on the generators $T$ and $S$ of $\tilde{\Gamma}$. We have $-\overline{M\tau} = M'(-\overline{\tau})$. From there on, the conditions for $T\psi$ to be a Maass cusp form can be checked straightforwardly.
\end{proof}
\begin{lemma}
For the Siegel theta function, we have that $T\circ \Theta_{S, \Delta}(z,\tau)$ is an automorphic form of weight $-1/2$ for $\overline{\tilde{\rho}}_L$, while for the Millson theta function, we have that $T\circ \Theta_{M, \Delta}(z, \tau)$ is an automorphic form of weight $1/2$ for $\overline{\tilde{\rho}}_L$.
\end{lemma}

\begin{remark}
The functions $T\circ\Theta_{S, \Delta}(z,\tau)$ and $T\circ\Theta_{M, \Delta}(z,\tau)$ are, in fact, the Siegel and Millson theta function corresponding to the lattice $L^-$ ($L$ regarded as a lattice in $(V,-q)$) of signature $(1,2)$. See, for example, \cite{swag} for the details regarding both theta functions for the signature $(1,2)$.
\end{remark}

Thus, we define the twisted Siegel theta lift for $\varphi \in U_{k,N}(s)$ of weight $0< k \in 4\mathbb{Z}+2$ by
\begin{align*}
I^*_{S, \Delta}(\tau, \varphi) := R^{(k+2)/4}_{-1/2, \tau}\bigg[T\circ \int\limits_{\Gamma\setminus\mathbb{H}}L^{k/2}_{k,z}\varphi(z)\Theta_{S, \Delta}(z,\tau)\frac{dxdy}{y^2}\bigg],
\end{align*}
and the twisted Millson theta lift for $\varphi \in U_{k,N}(s)$ of weight $0< k \in 4\mathbb{Z}$ by
\begin{align*}
I^*_{M, \Delta}(\tau, \varphi) := R^{k/4}_{1/2, \tau}\bigg[T\circ \int\limits_{\Gamma\setminus\mathbb{H}}L^{k/2}_{k,z}\varphi(z)\Theta_{M, \Delta}(z,\tau)\frac{dxdy}{y^2}\bigg].
\end{align*}
For weights $k < 0$, this generalizes as in Section \ref{lifts}.\\
Note that for any weight $k$ Maass cusp form $\psi(\tau)$, we have $R_{-k,\tau}[T\circ\psi(\tau)] = T\circ[L_{k,\tau}\psi(\tau)]$, as can be verified via a straightforward calculation on the Fourier series of both forms. This means that the Fourier expansion of the functions $I^*_{S, \Delta}(\tau, \varphi)$ and $I^*_{M, \Delta}(\tau, \varphi)$ can be derived from that of $I_{S, \Delta}(\tau, \varphi)$ and $I_{M, \Delta}(\tau, \varphi)$. We explain the procedure for the Siegel theta lift for non-negative weights, as the other cases are sufficiently similar. We have
\begin{align*}
&I^*_{S, \Delta}(\tau, \varphi) = R^{(k+2)/4}_{-1/2, \tau}\bigg[T\circ \int\limits_{\Gamma\setminus\mathbb{H}}L^{k/2}_{k,z}\varphi(z)\Theta_{S, \Delta}(z,\tau)\frac{dxdy}{y^2}\bigg]\\
&= T\circ \bigg[L^{(k+2)/4}_{1/2, \tau}\int\limits_{\Gamma\setminus\mathbb{H}}L^{k/2}_{k,z}\varphi(z)\Theta_{S, \Delta}(z,\tau)\frac{dxdy}{y^2}\bigg].
\end{align*}
According to the same arguments as in the proof of Theorem \ref{fünf}, the Fourier expansion of the expression in the brackets
is given by
\begin{align*}
\tiny{ \sum\limits_{h \in L'/L}\sum\limits_{\substack {D \in \mathbb{Z}+\textup{sign}(\Delta)q(h) \\D \neq 0}}\tilde{C}(h, D)(4\pi |D|v)^{\frac{k+1}{4}}W_{\textup{sign}(D)\frac{-(k+1)}{4}, \frac{s}{2}- \frac{1}{4}}(4\pi |D|v)e(nu)\mathfrak{e}_h,}
\end{align*}
where the coefficients $\tilde{C}(h, D)$ are given by
\begin{align*}
\tilde{C}(h, D) = (-1)^{\frac{k+2}{4}}(4\pi|D|)^{-\frac{k+2}{4}}\frac{\sqrt{\pi}}{2}|D|^{-1/2}Tr_{\Delta, \mathfrak{r}}(L_{k,z}^{-k/2}\varphi, h, D),
\end{align*}
for $D < 0$ and by
\begin{align*}
\tilde{C}(h, D) = (-1)^{\frac{k}{2}}(4\pi|D|)^{-\frac{k+2}{4}}2^{-\frac{k+2}{2}-1}\frac{\Gamma(s+k/2)}{\Gamma(s-k/2)}\cdot Tr_{\Delta, \mathfrak{r},k-2}(L_{k,z}\varphi, h, D),
\end{align*}
for $D > 0$.
\begin{remark}
Note that the differential operator in the cycle integral remains, due to the fact that $4 \nmid \frac{k}{2}$, so instead of regarding $\varphi$, we regard $L_{k,z}\varphi$ when running through the steps of the proof of Theorem \ref{fünf}.
\end{remark}
Now, we apply the operator $T$ to this. Executing this procedure in all possible cases, we end up with:
\begin{theorem}\label{sieben}
For $\varphi \in U_{k, N}(s)$ an even Maass cusp form with $0 \leq k \in 4\mathbb{Z}+2$, the Siegel theta lift $I^*_{S, \Delta}(\tau, \varphi)$ is given by
\begin{align*}
&I^*_{S, \Delta}(\tau, \varphi) = T\circ\tilde{\phi}_{\Delta}(\tau, \varphi)\\
&\tiny{=  \sum\limits_{h \in L'/L}\sum\limits_{\substack {D \in \mathbb{Z}-\textup{sign}(\Delta)q(h) \\D \neq 0}}\hat{C}(h,D)(4\pi |D|v)^{-\frac{k+1}{4}}W_{\textup{sign}(D)\frac{(k+1)}{4}, \frac{s}{2}- \frac{1}{4}}(4\pi |D|v)e(Du)\mathfrak{e}_h,}
\end{align*}
with
\begin{align*}
&\hat{C}(h,D) = \tilde{C}(h, -D)(4\pi |D|)^{\frac{k+1}{2}}\\
&= \begin{cases}
(-1)^{\frac{k+2}{4}}(4\pi|D|)^{\frac{k}{4}}\frac{\sqrt{\pi}}{2}|D|^{-1/2}Tr_{\Delta, \mathfrak{r}}(L_{k,z}^{-k/2}\varphi, h, -D) &\textup{ for } D>0,\\
 (-1)^{\frac{k+2}{4}}(4\pi|D|)^{\frac{k}{4}}2^{-\frac{k+2}{2}-1}\frac{\Gamma(s+k/2)}{\Gamma(s-k/2)}\cdot Tr_{\Delta, \mathfrak{r},k-2}(L_{k,z}\varphi, h,- D) &\textup{ for } D<0.
\end{cases}
\end{align*}
For $k < 0, k \in 4\mathbb{Z} +2$, we have
\begin{align*}
&I^*_{S, \Delta}(\tau, \varphi)\\
&\tiny{ =   \sum\limits_{h \in L'/L}\sum\limits_{\substack {D \in \mathbb{Z}-\textup{sign}(\Delta)q(h) \\D \neq 0}}\hat{C}(h,D)(4\pi |D|v)^{-\frac{k+1}{4}}W_{\textup{sign}(D)\frac{(k+1)}{4}, \frac{s}{2}- \frac{1}{4}}(4\pi |D|v)e(Du)\mathfrak{e}_h,}
\end{align*}
with
\begin{align*}
&\hat{C}(h,D)\\
&= \begin{cases}
(-1)^{\frac{k+2}{4}}(4\pi|D|)^{\frac{k}{4}}\frac{\sqrt{\pi}}{2}|D|^{-1/2}\frac{\Gamma\left(\frac{s}{2}-\frac{k}{4}\right)\Gamma\left(\frac{s}{2}\right)}{\Gamma\left(\frac{s}{2}+\frac{1}{2}\right)\Gamma\left(\frac{s}{2}+\frac{k+2}{4}\right)}\\\times Tr_{\Delta, \mathfrak{r}}(R_{k,z}^{-k/2}\varphi, h, -D) &\textup{ for } D>0,\\
 (-1)^{\frac{k+2}{4}}(4\pi|D|)^{\frac{k}{4}}2^{\frac{-k}{2}}\frac{\Gamma\left(\frac{s}{2}-\frac{k}{4}\right)\Gamma\left(\frac{s}{2}\right)}{\Gamma\left(\frac{s}{2}+\frac{1}{2}\right)\Gamma\left(\frac{s}{2}+\frac{k+2}{4}\right)}\\
\times Tr_{\Delta, \mathfrak{r},k-2}(R_{k,z}\varphi, h, -D) &\textup{ for } D<0.
\end{cases}
\end{align*}
\end{theorem}
\begin{theorem}\label{acht}
 For $\varphi \in U_{k, N}(s)$ an even Maass cusp form with $0 \leq k \in 4\mathbb{Z}$, the Millson theta lift $I^*_{M, \Delta}(\tau, \varphi)$ is given by
\begin{align*}
&I^*_{M, \Delta}(\tau, \varphi)\\
&\tiny{=  \sum\limits_{h \in L'/L}\sum\limits_{\substack {D \in \mathbb{Z}-\textup{sign}(\Delta)q(h) \\D \neq 0}}\hat{C}(h,D)(4\pi |D|v)^{-\frac{k+1}{4}}W_{\textup{sign}(D)\frac{(k+1)}{4}, \frac{s}{2}- \frac{1}{4}}(4\pi |D|v)e(Du)\mathfrak{e}_h,}
\end{align*}
with
\begin{align*}
&\hat{C}(h,D)\\
&= \begin{cases}
(-1)^{\frac{k}{4}}\frac{(4\pi|D|)^{\frac{k+2}{4}}}{4\sqrt{\pi}|D|} Tr^{\pm}_{\Delta, \mathfrak{r}}(L_{k,z}^{-k/2}\varphi, h, -D)  &\textup{ for } D>0,\\
 (-1)^{\frac{k}{4}} \frac{(4\pi|D|)^{\frac{k+2}{4}}}{8\pi D}2^{-\frac{k}{2}}\frac{\Gamma(s+k/2)}{\Gamma(s-k/2)}\cdot Tr_{\Delta, \mathfrak{r},k}(L_{k,z}\varphi, h, -D) &\textup{ for } D<0.
\end{cases}
\end{align*}
For $k < 0, k \in 4\mathbb{Z}$, we have
\begin{align*}
&I^*_{M, \Delta}(\tau, \varphi)\\
&\tiny{ =   \sum\limits_{h \in L'/L}\sum\limits_{\substack {D \in \mathbb{Z}-\textup{sign}(\Delta)q(h) \\D \neq 0}}\hat{C}(h,D)(4\pi |D|v)^{-\frac{k+1}{4}}W_{\textup{sign}(D)\frac{(k+1)}{4}, \frac{s}{2}- \frac{1}{4}}(4\pi |D|v)e(Du)\mathfrak{e}_h,}
\end{align*}
with
\begin{align*}
&\hat{C}(h,D)\\
&= \begin{cases}
(-1)^{\frac{k}{4}}\frac{(4\pi|D|)^{\frac{k+2}{4}}}{4\sqrt{\pi}|D|}\frac{\Gamma\left(\frac{s}{2}-\frac{k}{4}\right)\Gamma\left(\frac{s}{2}+\frac{1}{2}\right)}{\Gamma\left(\frac{s}{2}\right)\Gamma\left(\frac{s}{2}+\frac{k+2}{4}\right)} Tr^{\pm}_{\Delta, \mathfrak{r}}\left(R_{k,z}^{k/2}\varphi, h, -D\right) &\textup{ for } D>0,\\
 (-1)^{\frac{k}{4}}\frac{(4\pi|D|)^{\frac{k+2}{4}}2^{-\frac{k}{2}}}{8\pi D}\frac{\Gamma\left(\frac{s}{2}-\frac{k}{4}\right)\Gamma\left(\frac{s}{2}+\frac{1}{2}\right)}{\Gamma\left(\frac{s}{2}\right)\Gamma\left(\frac{s}{2}+\frac{k+2}{4}\right)}\cdot Tr_{\Delta, \mathfrak{r},k}(R_{k,z}\varphi, h, -D) &\textup{ for } D<0.
\end{cases}
\end{align*}
\end{theorem}

\section{The Shimura and Millson-Shimura Lift}\label{liftsec}
\subsection{Definition}
Let $F(\tau) \in U_{\frac{k+1}{2}, \tilde{\rho}_L}(s)$ be a Maass cusp form with Fourier series
\begin{align*}
 \sum\limits_{h \in L'/L}\big[\sum\limits_{\substack {n \in \mathbb{Z} + \textup{sign}(\Delta)q(h) \\n \neq 0}}c(h, n)(4\pi |n|v)^{-\frac{k+1}{4}}W_{\textup{sign}(n)\frac{k+1}{4}, s- \frac{1}{2}}(4\pi |n|v)e(nu)\big]\mathfrak{e}_h,
\end{align*}
and throughout this section, let $\mathcal{F} = \mathcal{F}_{\textup{SL}_2(\mathbb{Z})}$ be the fundamental domain of the full modular group.
\begin{definition}
Let $k \in 4\mathbb{Z}$, then for $k \geq 0$, we define the \textbf{Shimura lift} of $F(\tau)$ by
\begin{align*}
\Phi_{\Delta,}(F, z) := R_{0,z}^{k/2}\int\limits_{\mathcal{F}}\big\langle L_{\frac{k+1}{2}, \tau}^{k/4}F(\tau), \Theta_{S,\Delta}(z, \tau) \big\rangle v^{1/2}d\mu(\tau),
\end{align*}
and for $k < 0$ by
\begin{align*}
\Phi_{\Delta}(F, z) := L_{0, z}^{-k/2}\int\limits_{\mathcal{F}}\big\langle R_{\frac{k+1}{2}, \tau}^{-k/4}F(\tau), \Theta_{S,\Delta}(z, \tau) \big\rangle v^{1/2}d\mu(\tau).
\end{align*}
Now let $k \in 4\mathbb{Z} +2$. Then we define the \textbf{Millson-Shimura} lift for $k > 0$ by
\begin{align*}
\Psi_{\Delta}(F, z) := R_{0,z}^{k/2}\int\limits_{\mathcal{F}}\big\langle L_{\frac{k+1}{2}, \tau}^{(k+2)/4}F(\tau), \Theta_{M,\Delta}(z, \tau) \big\rangle v^{-1/2}d\mu(\tau)
\end{align*}
and for $k < 0$ by
\begin{align*}
\Psi_{\Delta}(F, z) := L_{0,z}^{-k/2}\int\limits_{\mathcal{F}}\big\langle R_{\frac{k+1}{2}, \tau}^{(-k-2)/4}F(\tau), \Theta_{M,\Delta}(z, \tau) \big\rangle v^{-1/2}d\mu(\tau).
\end{align*}
\end{definition}
It is clear that the functions $\Phi_{\Delta}(F, z)$ and $\Psi_{\Delta}(F, z)$, up to convergence of the integrals, transform like modular forms of weight $k$.
\begin{proposition}\label{ka}
If  $N$ is square-free and $F(\tau)$ lies in $U_{1/2, \tilde{\rho}_L}(\frac{s}{2}+ \frac{1}{4})$ and is non-harmonic $($i.e.\ $\lambda = \left(\frac{s}{2}+\frac{1}{4} - \frac{k+1}{4}\right)\left(1- \frac{s}{2}-\frac{1}{4} - \frac{k+1}{4}\right) \neq 0$$)$, then $\Phi_{\Delta}(F, z)$ lies in  $U_{0,N}(s)$.\\
If $N$ is square-free and $F(\tau)$ lies in $U_{-1/2, \tilde{\rho}_L}(\frac{s}{2}+ \frac{1}{4})$, then $\Psi_{\Delta}(F, z)$ lies in  $U_{0,N}(s)$.

\end{proposition}
\begin{proof}
Since one of the functions in the integral are of exponential decay towards the cusp of $\textup{SL}_2(\mathbb{Z})$ and the other at most grows like $\mathcal{O}(v^{1/2})$, the integral converges by the same arguments as in \ref{baum}.\\
For the Laplace property, we again use Lemma \ref{lapla} and Lemma 4.3 of \cite{MR1903920}.
We, furthermore, need to verify that $\Psi_{\Delta}(F, z)$ and $\Phi_{\Delta}(F,z)$ are indeed of exponential decay towards the cusps. This follows as in Proposition \ref{baum}.
\end{proof}

\begin{corollary}
If $F(\tau) \in U_{\frac{k+1}{2}, \tilde{\rho}_L}(\frac{s}{2} + \frac{1}{4})$ for $k \in 4\mathbb{Z}$ and does not become harmonic under $L_{\frac{k+1}{2},\tau}^{k/4}$ $($respectively $R_{\frac{k+1}{2},\tau}^{-k/4}$$)$ then the Shimura lift $\Phi_{\Delta}(F, z)$ lies in $U_{k,N}(s)$ and if  $F(\tau) \in U_{\frac{k+1}{2}, \tilde{\rho}_L}(\frac{s}{2} + \frac{1}{4})$ for $k \in 4\mathbb{Z}+2$, the Millson-Shimura lift $\Psi_{\Delta}(F, z)$ lies in $U_{k,N}(s)$.
\end{corollary}
\begin{proof}
This follows from Proposition \ref{ka} and the fact that $L_k$ and $R_k$ map Laplace eigenforms to Laplace eigenforms.%Lemma unlabeled
\end{proof}

\subsection{Fourier Coefficients}
We give a formula for the Fourier expansion of the Shimura and Millson-Shimura Lift. Our computations mainly follow the methods of \cite{MR0364106}, \cite{MR1244668}, \cite{hoe}, \cite{craw} and \cite{MR2726107}. From now on, we assume that $s \neq \ell/2$ with $\ell \equiv k \textup{ mod }2$. This implies that a function $F(\tau) \in U_{k,\tilde{\rho}_L}\left( \frac{s}{2}+\frac{1}{4}\right)$ does not become harmonic under the iterated lowering or raising operator.
\begin{theorem}\label{MShim}
Let $\varepsilon = 1$ when $\Delta > 0$ and $i$ when $\Delta < 0$.\\
Let  $F(\tau) \in U_{\frac{k+1}{2}, \tilde{\rho}_L}\left( \frac{s}{2}+\frac{1}{4}\right)$. Denote by $c(h, n)$ the $n$-th Fourier coefficient of $F(\tau)$. Then for $0 \leq k \in 4\mathbb{Z}$ the Shimura lift $\Phi_{\Delta}(F,z)$ is in $U_{k,N}(s)$ and has the Fourier expansion
\begin{align*}
&\Phi_{\Delta}(F,z) = (-1)^{k/4}2\varepsilon|\Delta|^{1/2} \frac{\Gamma\Big(\frac{s}{2}+\frac{k+2}{4}\Big)\Gamma\Big(\frac{s}{2}\Big)}{\Gamma\Big(\frac{s}{2}-\frac{k}{4} \Big) \Gamma\Big(\frac{s}{2}+\frac{1}{2}\Big)}\\
&\times \sum\limits_{n\neq 0}(4\pi ny)^{-k/2}a(n) W_{0, s-\frac{1}{2}}(4\pi |n|y)e^{2\pi i n x},
\end{align*}
where
\begin{align*}
a(n) =  \begin{cases}(-4\pi n)^{k/2}\sum\limits_{m\mid n}\frac{\begin{pmatrix}\frac{\Delta}{m}\end{pmatrix}}{m\big(4\pi\big|\Delta\frac{(n/m)^2}{4N}\big|\big)^{\frac{k}{4}}}c\bigg(\mathfrak{r}\frac{(n/m)}{2N}, |\Delta|\frac{(n/m)^2}{4N}\bigg), \textup{ for } n > 0,\\
(-4\pi n)^{k/2}\sum\limits_{m\mid n}\frac{D_{s,k}\begin{pmatrix}\frac{\Delta}{m}\end{pmatrix}}{m\big(4\pi\big|\Delta\frac{(n/m)^2}{4N}\big|\big)^{\frac{k}{4}}}c\bigg(\mathfrak{r}\frac{(n/m)}{2N}, |\Delta|\frac{(n/m)^2}{4N}\bigg), \textup{ for } n < 0
\end{cases}
\end{align*}
with $D_{s,k} = \frac{\Gamma\big(s+\frac{k}{2}\big)}{\Gamma\big(s - \frac{k}{2}\big)}$. \\
For $k < 0, k \in 4\mathbb{Z}$, the Shimura lift $\Phi_{\Delta}(F,z)$ in $U_{k,N}(s)$ and has the Fourier expansion
\begin{align*}
&\Phi_{\Delta}(F,z) = (-1)^{k/4}2\varepsilon|\Delta|^{1/2} \sum\limits_{n\neq 0}(4\pi ny)^{k/2}a(n) W_{0, s-\frac{1}{2}}(4\pi |n|y)e^{2\pi i n x},
\end{align*}
where
\begin{align*}
a(n) =  \begin{cases}(-4\pi n)^{\frac{k}{2}}\sum\limits_{m\mid n}\frac{\tilde{D}_{s,k}\begin{pmatrix}\frac{\Delta}{m}\end{pmatrix}}{m\big(4\pi\big|\Delta\frac{(n/m)^2}{4N}\big|\big)^{\frac{k}{4}}}c\bigg(\mathfrak{r}\frac{(n/m)}{2N}, |\Delta|\frac{(n/m)^2}{4N}\bigg), \textup{ for } n > 0,\\
(-4\pi n)^{\frac{k}{2}}\sum\limits_{m\mid n}\frac{\begin{pmatrix}\frac{\Delta}{m}\end{pmatrix}}{m \big(4\pi\big|\Delta\frac{(n/m)^2}{4N}\big|\big)^{\frac{k}{4}}}c\bigg(\mathfrak{r}\frac{(n/m)}{2N}, |\Delta|\frac{(n/m)^2}{4N}\bigg), \textup{ for } n < 0
\end{cases}
\end{align*}
with $\tilde{D}_{k,s} =  \frac{\Gamma\big(s-\frac{k}{2}\big)}{\Gamma\big(s +\frac{k}{2}\big)}$.\\
 For $0 < k \in 4\mathbb{Z} + 2$, the Millson-Shimura lift $\Psi_{\Delta}(F, z)$ is  in $U_{k,N}(s)$ and has the Fourier expansion
\begin{align*}
&\Psi_{\Delta}(F,z) = (-1)^{\frac{k+2}{4}}\frac{2\varepsilon|\Delta|}{i\sqrt{N}} \frac{\Gamma\Big(\frac{s}{2}+\frac{k+2}{4}\Big)\Gamma\Big(\frac{s}{2} + \frac{1}{2}\Big)}{\Gamma\Big(\frac{s}{2}-\frac{k}{4} \Big)\Gamma\big( \frac{s}{2}\big)}\\
&\times \sum\limits_{n\neq 0}(4\pi ny)^{-k/2}a(n) W_{0, s-\frac{1}{2}}(4\pi |n|y)e^{2\pi i n x}, %k > 0
\end{align*}
where
\begin{align*}
a(n) =  \begin{cases}(-4\pi n)^{\frac{k}{2}}\sum\limits_{m\mid n}\frac{\begin{pmatrix}\frac{\Delta}{m}\end{pmatrix}n}{m^2\big(4\pi\big|\Delta\frac{(n/m)^2}{4N}\big|\big)^{\frac{k+2}{4}}}c\bigg(\mathfrak{r}\frac{(n/m)}{2N}, |\Delta|\frac{(n/m)^2}{4N}\bigg), \textup{ for } n > 0,\\
(-4\pi n)^{\frac{k}{2}}\sum\limits_{m\mid n}\frac{D_{s,k}\begin{pmatrix}\frac{\Delta}{m}\end{pmatrix}n}{m^2\big(4\pi\big|\Delta\frac{(n/m)^2}{4N}\big|\big)^{\frac{k+2}{4}}}c\bigg(\mathfrak{r}\frac{(n/m)}{2N}, |\Delta|\frac{(n/m)^2}{4N}\bigg), \textup{ for } n < 0.
\end{cases}
\end{align*}
And for $k < 0, k \in 4\mathbb{Z} +2$, the Millson-Shimura lift $\Psi_{\Delta}(F, z)$ is  in $U_{k,N}(s)$ and has the Fourier expansion
\begin{align*}
&\Psi_{\Delta}(F, z) = (-1)^{\frac{k+2}{4}}\frac{2\varepsilon|\Delta|}{i\sqrt{N}} \sum\limits_{n\neq 0}(4\pi ny)^{k/2}a(n) W_{0, s-\frac{1}{2}}(4\pi |n|y)e^{2\pi i n x}, %k > 0
\end{align*}
with coefficients
\begin{align*}
a(n) =  \begin{cases}(-4\pi n)^{\frac{k}{2}}\sum\limits_{m\mid n}\frac{\tilde{D}_{s,k}\begin{pmatrix}\frac{\Delta}{m}\end{pmatrix}n}{m^2\big(4\pi\big|\Delta\frac{(n/m)^2}{4N}\big|\big)^{\frac{k+2}{4}}}c\bigg(\mathfrak{r}\frac{(n/m)}{2N}, |\Delta|\frac{(n/m)^2}{4N}\bigg), \textup{ for } n > 0,\\
(-4\pi n)^{\frac{k}{2}}\sum\limits_{m\mid n}\frac{\begin{pmatrix}\frac{\Delta}{m}\end{pmatrix}n}{m^2\big(4\pi\big|\Delta\frac{(n/m)^2}{4N}\big|\big)^{\frac{k+2}{4}}}c\bigg(\mathfrak{r}\frac{(n/m)}{2N}, |\Delta|\frac{(n/m)^2}{4N}\bigg), \textup{ for } n < 0.
\end{cases}
\end{align*}
\end{theorem}
\begin{proof}
The proof mainly follows the idea of Niwa, \cite{MR0364106}. We give the proof for the Millson Theta function as the other case follows by the same methods. Lets first assume that $k = 0$. We regard the Mellin transform of the Millson-Shimura lift, given by
\begin{align*}
&\Omega(r) = \int\limits_{0}^{\infty}\Psi_{\Delta}(F, i\eta)\eta^{r-1}d\eta\\
& =  \int\limits_{0}^{\infty}\int\limits_{\mathcal{F}}\langle F(\tau), \Theta_{M, \Delta}(i\eta, \tau) \rangle v^{-1/2}d\mu(\tau)\eta^{r-1}d\eta =  \int\limits_{0}^{\infty}\int\limits_{\mathcal{F}}\bigg\langle F(\tau),\sum\limits_{h\in L'/L}\\&\bigg[ \frac{iN\eta^2}{|\Delta|}\sum\limits_{c,d \in \mathbb{Z}}(c\tau + d)e^{-\frac{\pi N\eta^2|c\tau +d|^2}{|\Delta|v}}\xi_h(\tau, \mu(i\eta), d, -c) \bigg]\bigg\rangle v^{-1/2}d\mu(\tau)\eta^{r-1}d\eta,
\end{align*}
where the last identity is due to Proposition \ref{thetah}. By regarding Lemma \ref{ll}, we obtain that this expression is given by
\begin{align*}
&\int\limits_{0}^{\infty}\int\limits_{\mathcal{F}}\bigg\langle F(\tau), \sum\limits_{h \in L'/L}\bigg[ \frac{iN\eta^2}{|\Delta|}
\sum\limits_{c,d \in \mathbb{Z}}(c\tau + d)e^{-\frac{\pi N\eta^2|c\tau +d|^2}{|\Delta|v}} \sum\limits_{\lambda \in K + \mathfrak{r}h}\\&\sum\limits_{\substack{t(\Delta)\\ q(\lambda +c\ell'+t\ell) \equiv \Delta q(h)(\Delta)}}\chi_{\Delta}(\lambda +c\ell'+t\ell)e\bigg(\frac{-dt}{|\Delta|}\bigg) e\bigg(\frac{\lambda ^2\tau}{2|\Delta|} \bigg)\bigg]\bigg\rangle v^{-1/2}d\mu(\tau)\eta^{r-1}d\eta.
\end{align*}
Note especially that the theta function $\xi_h(\tau, \mu\big(z), d, -c)$ only depends on $\tau$ if $z$ is on the imaginary axis, i.e $\mu(i\eta) = 0$. Following \cite{MR1244668}, we set
\begin{align*}
\beta := \frac{\sqrt{\pi N} \eta|c\tau + d|}{{\sqrt{v|\Delta|}}}
\hspace{2mm}
\textup{ and }
\eta = \frac{{\sqrt{v|\Delta|}\beta}}{\sqrt{\pi N}|c\tau + d|}.
\end{align*}
Thus, our expression is given by
\begin{align*}
&\small{\frac{N^{\frac{-r}{2}}}{\pi^{\frac{(r+2)}{2}}}|\Delta|^{\frac{r}{2}}\int\limits_0^{\infty} e^{-\beta^2}\beta^{r+1}d\beta}\\
&\small{\times \int\limits_{\mathcal{F}}\bigg\langle F(\tau), \sum\limits_{h \in L'/L}\sum\limits_{c,d \in \mathbb{Z}}i(c\tau + d)\bigg(\frac{v}{|c\tau+d|^2} \bigg)^{\frac{r+2}{2}}}\\
&\small{\sum\limits_{\lambda \in K + \mathfrak{r}h}\sum\limits_{\substack{t(\Delta)\\ q(\lambda +c \ell'+t\ell) \equiv \Delta q(h)(\Delta)}}\chi_{\Delta}(\lambda +c\ell'+t\ell)e(-dt/|\Delta|)
 e\bigg(\frac{\lambda ^2}{2}\frac{\tau}{|\Delta|} \bigg) \bigg]\bigg\rangle v^{-1/2}d\mu(\tau)}\\
&\small{ = \frac{1}{2}\frac{N^{\frac{-r}{2}}}{\pi^{\frac{(r+2)}{2}}}|\Delta|^{\frac{r}{2}}\Gamma\big(\frac{r+2}{2}\big)\int\limits_{\mathcal{F}}\bigg\langle F(\tau), \sum\limits_{h \in L'/L}\sum\limits_{c,d \in \mathbb{Z}}i(c\tau + d)\bigg(\frac{v}{|c\tau+d|^2} \bigg)^{\frac{r+2}{2}} }\\
 &\small{\sum\limits_{\lambda \in K + \mathfrak{r}h}\sum\limits_{\substack{t(\Delta)\\ q(\lambda  +c\ell'+t\ell) \equiv \Delta q(h)(\Delta)}}\chi_{\Delta}(\lambda +c\ell'+t\ell)e(-dt/|\Delta|)
 e\bigg(\frac{\lambda ^2}{2}\frac{\tau}{|\Delta|} \bigg) \bigg]\bigg\rangle v^{-1/2}d\mu(\tau).}
\end{align*}
The sum over all $c, d \in \mathbb{Z}$ can be rewritten as
\begin{align*}
&\sum\limits_{c,d \in \mathbb{Z}}i(c\tau + d)\bigg(\frac{v}{|c\tau+d|^2} \bigg)^{\frac{r+2}{2}}\xi_h(\tau, \mu, d, -c)\\
&= \sum\limits_{n\geq 1}n^{-(r+1)}\sum\limits_{(d, c) =1}i(c\tau + d)\bigg(\frac{v}{|c\tau+d|^2} \bigg)^{\frac{r+2}{2}}\xi_h(\tau, \mu, nd, -nc).
\end{align*}
Note that for $n = 0$ the function vanishes. Two matrices $\begin{psmallmatrix}a&b\\c&d \end{psmallmatrix}$ and $\begin{psmallmatrix}a'&b'\\c'&d' \end{psmallmatrix}$ are equal in $\Gamma_{\infty}\setminus \Gamma$ if and only if their bottom row is the same, i.e.\ $d = d'$ and $c=c'$. Furthermore, any tuple $(c,d) \in \mathbb{Z}^2$ with $\textup{gcd}(c,d) = 1$ is the bottom row of an element in $\Gamma$ and, trivially, all elements of $\Gamma$ fulfill this condition. Due to Proposition \ref{poisson}, we know that for $(M, \phi) = \big(\begin{psmallmatrix}a&b\\c&d\end{psmallmatrix}, \phi\big) \in \tilde{\Gamma}$, we have
\begin{align*}
\Xi(\tau, \mu, nd  ,  -nc )= \phi(\tau)^{-1}\tilde{\rho}_K(M, \phi)^{-1}\Xi(M\tau, \mu, n, 0).
\end{align*}
Furthermore, we have that $\textup{Im}(M\tau) = \frac{\textup{Im}(\tau)}{|c\tau + d|^2}$ and $\phi_M(\tau)^2 = c\tau + d$. Taking all this into account, we can rewrite the sum over all $(c, d) =1$ as a sum over $\tilde{\Gamma}_{\infty}\setminus \tilde{\Gamma}$. Note that one has to consider the contribution of $(M, \phi_M)$ and $(M, -\phi_M)$, thus we obtain the factor $1/2$ in the following. We  get
\begin{align*}
&\small{\sum\limits_{h \in L'/L}\sum\limits_{c,d \in \mathbb{Z}}i(c\tau + d)\bigg(\frac{v}{|c\tau+d|^2} \bigg)^{\frac{r+2}{2}}\xi_h(\tau, 0, d, -c)} \\
&\small{=\sum\limits_{n\geq 1}n^{-(r+1)}\sum\limits_{(c,d)=1}i(c\tau + d)\bigg(\frac{v}{|c\tau+d|^2} \bigg)^{\frac{r+2}{2}}\phi(\tau)^{-1}\tilde{\rho}_K(M, \phi)^{-1}\Xi(M\tau, 0, n, 0)}\\
&\small{ =\sum\limits_{n\geq 1}n^{-(r+1)}\sum\limits_{M' \in \tilde{\Gamma}_{\infty}\setminus \tilde{\Gamma}}\frac{i}{2}\big[v^{\frac{r+2}{2}} \cdot \Xi(\tau, 0, n, 0)\big]\Big|_{-1/2, \tilde{\rho}_L}M'}.
\end{align*}
For the remaining integral in our Mellin transform this means that it is given by
\begin{align*}
\int\limits_{\mathcal{F}}\Big\langle F(\tau),  \sum\limits_{n\geq 1}n^{-(r+1)}\sum\limits_{M' \in \tilde{\Gamma}_{\infty}\setminus \tilde{\Gamma}}\frac{i}{2}\big[v^{\frac{r+2}{2}} \cdot \Xi(\tau, 0, n, 0)\big]\Big|_{-1/2, \tilde{\rho}_L}M'\Big\rangle v^{-1/2}d\mu(\tau).
\end{align*}
Recall that the function in the integral as a whole is automorphic of weight $0$, so this expression can be written as a sum of integrals
\begin{align*}
&-2\int\limits_{\mathcal{F}} \sum\limits_{n\geq 1}in^{-(r+1)}\textup{Im}(\tau)^{\frac{r+1}{2}}\big\langle F(\tau), \Xi(\tau, 0, n, 0)\big\rangle d\mu(\tau)\\
&- \int\limits_{\mathcal{F}} \sum\limits_{n\geq 1}in^{-(r+1)}\sum\limits_{\substack{M =\begin{psmallmatrix}a&b\\c&d\end{psmallmatrix} \in \Gamma_{\infty}\setminus \Gamma\\ c\neq 0}} \textup{Im}(M\tau)^{\frac{r+1}{2}}\langle F(\tau), \Xi(M\tau, 0, n, 0)\big\rangle d\mu(\tau).
\end{align*}
Here, we have split the sum over all matrices in $\tilde{\Gamma}_{\infty}\setminus \tilde{\Gamma}$ into that over all $M = \begin{psmallmatrix}a&b\\c&d\end{psmallmatrix}$ with $c \neq 0$ and that over $\pm\begin{psmallmatrix}1&0\\0&1 \end{psmallmatrix}$ and rewritten it as a sum over $\Gamma_{\infty}\setminus \Gamma$, which -- again by considering $(M, \phi)$ and $(M, -\phi)$ -- gives the factor $2$ in the second integral and the factor $4$ in the first, since one has to consider the contribution of both, $\pm M$. Due to the cuspidal properties of $F(\tau)$, both integrals converge absolutely. This now allows us to apply the usual unfolding argument (Rankin-Selberg method, see \cite{MR1625724} or \cite{MR1903920}). We obtain that the second integral is equal to
\begin{align*}
-2 \int\limits_{\mathcal{G}} \sum\limits_{n\geq 1}n^{-(r+1)}i \textup{Im}(\tau)^{\frac{r+1}{2}}\langle F(\tau), \Xi(\tau, 0, n, 0)\big\rangle d\mu(\tau),
\end{align*}
where $\mathcal{G} = \{\tau \in \mathbb{H} \mid |\textup{Re}(\tau)| \leq 1/2, |\tau| \leq 1\}$ is a fundamental domain for the action of $\Gamma_{\infty}$ on $\mathbb{H}-\bigcup\limits_{M \in \Gamma_{\infty}}M\mathcal{F}$.  %mehr dazu schreiben!
Combining both integrals, we obtain that their sum is equal to
\begin{align*}
-2\int\limits_0^{\infty}\int\limits_{-1/2}^{1/2} \sum\limits_{n\geq 1}n^{-(r+1)}i \textup{Im}(\tau)^{\frac{r-3}{2}}\langle F(\tau), \Xi(\tau, 0, n, 0)\big\rangle dudv.
\end{align*}
Now, we insert the Fourier expansion of both forms, which  yields to
\small
\begin{align} \label{frr}
&\tiny{= -2i\varepsilon|\Delta|^{1/2}\int\limits_0^{\infty}\int\limits_{-1/2}^{1/2} \sum\limits_{n\geq 1}\frac{\begin{psmallmatrix}\frac{\Delta}{n}\end{psmallmatrix}}{n^{-(r+1)}}\sum\limits_{h \in K'/K}\sum\limits_{\substack{\lambda \in K + \mathfrak{r}h \\q(\lambda) = \Delta q(h)(\Delta)}}}\\
&\tiny{\sum\limits_{\substack {m \in \mathbb{Z} + \textup{sign}(\Delta)q(h) \\m \neq 0}}c(h, m)(4\pi |m|v)^{\frac{1}{4}}W_{\textup{sign}(m)\frac{-1}{4}, s-\frac{1}{2}}(4\pi |m|v)e(mu)\overline{e\bigg(\frac{\lambda^2\tau}{2|\Delta|}\bigg)}v^{\frac{s-3}{2}}dudv} \notag \\
&\tiny{= - 2i\varepsilon|\Delta|^{1/2}L(\chi'_{\Delta}, r+1)\int\limits_0^{\infty}\int\limits_0^1 \sum\limits_{h \in K'/K}\sum\limits_{\substack{\lambda \in K + \mathfrak{r}h \\q(\lambda) = \Delta q(h)(\Delta)}}} \notag \\
&\tiny{\sum\limits_{\substack {m \in \mathbb{Z} + \textup{sign}(\Delta)q(h) \\m \neq 0}} c(h, m)(4\pi |m|v)^{\frac{1}{4}}W_{\textup{sign}(m)\frac{-1}{4}, s-\frac{1}{2}}(4\pi |m|v)e(mu)\overline{e\bigg(\frac{\lambda^2\tau}{2|\Delta|}\bigg)}v^{\frac{r-3}{2}}dudv,}  \notag %chi beibehalten
\end{align}
\normalsize
where $\chi'_{\Delta}$ denotes the Dirichlet character $\begin{psmallmatrix} \frac{\Delta}{\cdot}\end{psmallmatrix}$ and $L(\chi'_{\Delta}, r)$ its corresponding $L$-function. We calculate the integral over $u$. By definition of the sums, the term $m - \frac{q(\lambda)}{|\Delta|}$ is always integral, thus, the integral over $u$ vanishes unless $m = \frac{q(\lambda)}{|\Delta|}$. We obtain
\begin{align*}
&-2i\varepsilon|\Delta|^{1/2}L(\chi'_{\Delta}, r+1)\int\limits_0^{\infty} \sum\limits_{h \in K'/K}\sum\limits_{\substack{\lambda \in K + \mathfrak{r}h \\q(\lambda) = \Delta q(h)(\Delta)}}c\bigg(h, \frac{\lambda^2}{2|\Delta|}\bigg)\bigg(4\pi \bigg|\frac{\lambda^2}{2|\Delta|}\bigg|v\bigg)^{\frac{1}{4}}\\
&\times W_{\textup{sign}\big(\frac{q(\lambda)}{|\Delta|}\big)\frac{-1}{4}, s-\frac{1}{2}}\bigg(4\pi \bigg|\frac{\lambda^2}{2|\Delta|}\bigg|v\bigg)e\bigg(\frac{\lambda^2iv}{2|\Delta|}\bigg)v^{\frac{r-3}{2}}dv.
\end{align*}
We have $\Delta = \mathfrak{r}^2 (4N)$, where $4N$ is also the level of $K$. Thus, for $\lambda \in K', h\in K'/K$ we have that it follows from $\lambda \equiv \mathfrak{r}h (K)$ and $q(\lambda) = \Delta q(h) (\Delta)$ that $\lambda = \Delta\lambda'$ for $\lambda' \in K'$ and $\mathfrak{r}\lambda' \equiv h(K)$. The expression above then equals
\begin{align*}
&-2i\varepsilon|\Delta|^{1/2}L(\chi'_{\Delta}, r+1)\int\limits_0^{\infty}\sum\limits_{\lambda \in K'}c(\mathfrak{r}\lambda, |\Delta|q(\lambda))(4\pi ||\Delta|q(\lambda)|v)^{\frac{1}{4}}\\
&\times W_{\textup{sign}(|\Delta|q(\lambda))\frac{-1}{4}, s-\frac{1}{2}}(4\pi ||\Delta|q(\lambda)|v)e\big(|\Delta|q(\lambda)iv\big)v^{\frac{r-3}{2}}dv.
\end{align*}
Any sum over $\lambda \in K'$ can be understood as a sum over $k \in \mathbb{Z}$ with $q(k) = \frac{k^2}{4N}$. Note as well that $K$ is positive definite. We end up with
\begin{align*}
&-2i\varepsilon|\Delta|^{1/2}L(\chi'_{\Delta}, r+1)\int\limits_0^{\infty}\sum\limits_{k \in \mathbb{Z}}c\bigg(\mathfrak{r}\frac{k}{2N}, |\Delta|\frac{k^2}{4N}\bigg)\bigg(4\pi \bigg||\Delta|\frac{k^2}{4N}\bigg|v\bigg)^{\frac{1}{4}}\\
&\times W_{-\frac{1}{4}, s-\frac{1}{2}}\bigg(4\pi \bigg||\Delta|\frac{k^2}{4N}\bigg|v\bigg)e\bigg(\frac{|\Delta|k^2iv}{4N}\bigg)v^{\frac{r-3}{2}}dv.
\end{align*}
Taking into account the asymptotics of the Whittaker function and the $e(\cdot)$ term, by a standard argument involving dominated convergence, we can exchange summation and integration. By substituting $\frac{|\Delta|k^2}{4N}v$ with $v$, it then remains to compute
\begin{align*}
&-4i\varepsilon|\Delta|^{1/2}L(\chi'_{\Delta}, r+1) \sum\limits_{k\geq 1}\bigg(\frac{|\Delta|}{4N}\bigg)^{\frac{-r+1}{2}}\frac{c\bigg(\mathfrak{r}\frac{k}{2N}, |\Delta|\frac{k^2}{4N}\bigg)}{k^{r-1}}\\
&\times \int\limits_0^{\infty}(4\pi)^{1/4}W_{-\frac{1}{4}, s- \frac{1}{2}}(4\pi v)e^{-2\pi v}v^{\frac{r}{2}-\frac{1}{4}}\frac{dv}{v}.
\end{align*}
According to \cite{MR2360010}, p. 823, the remaining integral is given by
\begin{align*}
&\int\limits_0^{\infty}(4\pi)^{1/4}W_{-\frac{1}{4}, s- \frac{1}{2}}(4\pi v)e^{-2\pi v}v^{\frac{r}{2}-\frac{1}{4}}\frac{dv}{v}\\
&=(4\pi)^{\frac{-r+1}{2}}\frac{\Gamma\Big(\frac{r}{2}+\frac{1}{4}-\big(s-\frac{1}{2}\big)\Big)\Gamma\Big(\frac{r}{2}+\frac{1}{4}+\big(s-\frac{1}{2}\big)\Big)}{\Gamma\Big(\frac{r+2}{2}\Big)}.
\end{align*}
We set
\begin{align*}
a(n) = \sum\limits_{m\mid n}m^{-2}n\begin{pmatrix}\frac{\Delta}{m}\end{pmatrix}c\bigg(\mathfrak{r}\frac{(n/m)}{2N}, |\Delta|\frac{(n/m)^2}{4N}\bigg),
\end{align*}
with $n, m > 0$, so that $\Omega(r)$ can be written in the following way:
\begin{align*}
\Omega(r) = \frac{2\varepsilon\pi^{-r -\frac{1}{2}}}{i}\bigg(\frac{|\Delta|}{\sqrt{N}}\bigg)\Gamma\Big(\frac{r}{2}+\frac{1}{4}-\big(s-\frac{1}{2}\big)\Big)\Gamma\Big(\frac{r}{2}+\frac{1}{4}+\big(s-\frac{1}{2}\big)\Big)\sum\limits_{n\geq 1}\frac{a(n)}{n^r}.
\end{align*}
Comparing this to Proposition \ref{mellin}, we see that the expression above is the Mellin transform of an even Maass cusp form of weight $0$ with spectral point $2s-1/2$ and Fourier expansion given by
\begin{align*}
\Psi_{\Delta}(F,z) =  \frac{2\varepsilon|\Delta|}{i\sqrt{N}}\sum\limits_{n\neq 0}a(n) W_{0, 2s-1}(4\pi |n|y)e^{2\pi i n x},
\end{align*}
where, according to the definition of an even Maass cusp form, we have set $a(-n) = a(n)$.\\
The case of general weight follows with Lemma \ref{der} in the following steps (here exemplary for the $0 \leq k \in 4\mathbb{Z} +2$ case):\\
\begin{itemize}
\item compute the Fourier expansion of $L_{\frac{k+1}{2},\tau}^{(k+2)/4}F(\tau)$ with an iterated application of Lemma \ref{der},
\item compute the Fourier expansion of $\Psi_{\Delta}\left(L_{\frac{k+1}{2},\tau}^{(k+2)/4}F(\tau), z\right)$ with Theorem \ref{MShim},
\item compute the Fourier expansion of $R_{0,z}^{k/2}\Psi_{\Delta}\left(L_{\frac{k+1}{2},\tau}^{(k+2)/4}F(\tau), z\right)$ with another iterated application of Lemma \ref{der}.
\end{itemize}
Executing these steps (with the differential operators adjusted accordingly) gives us the statement.
\end{proof}

\subsection{The Shimura and Millson-Shimura Theta Lift for Switched Signature}
We will additionally give the formula for the Shimura and Millson-Shimura Lift for switched signature, which we need in order to obtain the full statement of Theorem \ref{KSmain}. However, we will not go into detail, since the changes to the proof of Theorem \ref{MShim} are minimal. For a detailed description on how to obtain the following corollary, regard \cite{str}. 
\begin{corollary}\label{cs}
Let $\varepsilon = 1$ when $\Delta > 0$ and $i$ when $\Delta < 0$.\\
Let  $F(\tau) \in U_{\frac{k+1}{2}, \overline{\tilde{\rho}}_L}\left( \frac{s}{2}+\frac{1}{4}\right)$. Denote by $c(h, n)$ the $n$-th Fourier coefficient of $F(\tau)$. Then for $0 \leq k \in 4\mathbb{Z}+2$ the Shimura lift $\Phi^*_{\Delta}(F,z)$ is in $U_{k,N}(s)$ and has the Fourier expansion
\begin{align*}
&\Phi^*_{\Delta}(F, z) = (-1)^{\frac{k+2}{4}} 2\sqrt{\pi} \varepsilon\bigg(\frac{|\Delta|^2}{N}\bigg)^{1/2}\sum\limits_{n\neq 0}(4\pi ny)^{k/2}a(n) W_{0, s-\frac{1}{2}}(4\pi |n|y)e^{2\pi i n x},
\end{align*}
where
\small
\begin{align*}
a(n) =  \begin{cases}(-4\pi n)^{k/2}\sum\limits_{m\mid n}\frac{n\begin{pmatrix}\frac{\Delta}{m}\end{pmatrix}}{m^2\big(4\pi\big|\Delta\frac{(n/m)^2}{4N}\big|\big)^{\frac{k+2}{4}}}c\bigg(\mathfrak{r}\frac{(n/m)}{2N}, -|\Delta|\frac{(n/m)^2}{4N}\bigg), \textup{ for } n > 0,\\
(-4\pi n)^{k/2}\sum\limits_{m\mid n}\frac{D_{s,k} n \begin{pmatrix}\frac{\Delta}{m}\end{pmatrix}}{m^2\big(4\pi\big|\Delta\frac{(n/m)^2}{4N}\big|\big)^{\frac{k+2}{4}}}c\bigg(\mathfrak{r}\frac{(n/m)}{2N}, -|\Delta|\frac{(n/m)^2}{4N}\bigg), \textup{ for } n < 0.
\end{cases}
\end{align*}
For $k < 0, k \in 4\mathbb{Z} +2$, the Shimura lift $\Phi^*_{\Delta}(F,z)$ is in $U_{k,N}(s)$ and has the Fourier expansion
\begin{align*}
&\Phi^*_{\Delta}(F,z) = (-1)^{\frac{k+2}{4}}2\sqrt{\pi} \varepsilon\bigg(\frac{|\Delta|^2}{N}\bigg)^{1/2}\frac{\Gamma\Big(\frac{s}{2}-\frac{k}{4}\Big)\Gamma\Big(\frac{s}{2}\Big)}{\Gamma\Big(\frac{s}{2}+\frac{k+2}{4} \Big)\Gamma\big( \frac{s}{2}+\frac{1}{2}\big)}\\
&\times\sum\limits_{n\neq 0}(4\pi ny)^{k/2}a(n) W_{0, s-\frac{1}{2}}(4\pi |n|y)e^{2\pi i n x},
\end{align*}
where
\small
\begin{align*}
a(n) =  \begin{cases}(-4\pi n)^{\frac{k}{2}}\sum\limits_{m\mid n}\frac{\tilde{D}_{s,k}n\begin{pmatrix}\frac{\Delta}{m}\end{pmatrix}}{m^2\big(4\pi\big|\Delta\frac{(n/m)^2}{4N}\big|\big)^{\frac{k+2}{4}}}c\bigg(\mathfrak{r}\frac{(n/m)}{2N}, -|\Delta|\frac{(n/m)^2}{4N}\bigg), \textup{ for } n > 0,\\
(-4\pi n)^{\frac{k}{2}}\sum\limits_{m\mid n}\frac{n\begin{pmatrix}\frac{\Delta}{m}\end{pmatrix}}{m^2\big(4\pi\big|\Delta\frac{(n/m)^2}{4N}\big|\big)^{\frac{k+2}{4}}}c\bigg(\mathfrak{r}\frac{(n/m)}{2N}, -|\Delta|\frac{(n/m)^2}{4N}\bigg), \textup{ for } n < 0.
\end{cases}
\end{align*}
 For $0< k \in 4\mathbb{Z}$, the Millson-Shimura lift $\Psi^*_{\Delta}(F,z)$ is in $U_{k,N}(s)$ and has the Fourier expansion
\begin{align*}
&\Psi^*_{\Delta}(F,z) = (-1)^{k/4} 2\frac{\varepsilon |\Delta|^{1/2}}{\sqrt{\pi} i} \sum\limits_{n\neq 0}(4\pi ny)^{-k/2}a(n) W_{0, s-\frac{1}{2}}(4\pi |n|y)e^{2\pi i n x}, %k > 0
\end{align*}
where
\begin{align*}
a(n) =  \begin{cases}(-4\pi n)^{\frac{k}{2}}\sum\limits_{m\mid n}\frac{\begin{pmatrix}\frac{\Delta}{m}\end{pmatrix}}{m\big(4\pi\big|\Delta\frac{(n/m)^2}{4N}\big|\big)^{\frac{k}{4}}}c\bigg(\mathfrak{r}\frac{(n/m)}{2N}, -|\Delta|\frac{(n/m)^2}{4N}\bigg), \textup{ for } n > 0,\\
(-4\pi n)^{\frac{k}{2}}\sum\limits_{m\mid n}\frac{D_{s,k}\begin{pmatrix}\frac{\Delta}{m}\end{pmatrix}}{m\big(4\pi\big|\Delta\frac{(n/m)^2}{4N}\big|\big)^{\frac{k}{4}}}c\bigg(\mathfrak{r}\frac{(n/m)}{2N}, -|\Delta|\frac{(n/m)^2}{4N}\bigg), \textup{ for } n < 0.
\end{cases}
\end{align*}
and, finally, for $k < 0, k \in 4\mathbb{Z}$, the Millson-Shimura lift $\Psi^*_{\Delta}(F,z)$ is in $U_{k,N}(s)$ and has the Fourier expansion
\begin{align*}
&\Psi^*_{\Delta}(F,z) = (-1)^{k/4} 2\frac{\varepsilon |\Delta|^{1/2}}{\sqrt{\pi} i}\frac{\Gamma\Big(\frac{s}{2}-\frac{k}{4}\Big)\Gamma\Big(\frac{s}{2}+\frac{1}{2}\Big)}{\Gamma\Big(\frac{s}{2}+\frac{k+2}{4} \Big)\Gamma\big( \frac{s}{2}\big)}\\
&\times \sum\limits_{n\neq 0}(4\pi ny)^{k/2}a(n) W_{0, s-\frac{1}{2}}(4\pi |n|y)e^{2\pi i n x}, %k > 0
\end{align*}
with coefficients
\begin{align*}
a(n) =  \begin{cases}(-4\pi n)^{\frac{k}{2}}\sum\limits_{m\mid n}\frac{\tilde{D}_{s,k}\begin{pmatrix}\frac{\Delta}{m}\end{pmatrix}}{m\big(4\pi\big|\Delta\frac{(n/m)^2}{4N}\big|\big)^{\frac{k}{4}}}c\bigg(\mathfrak{r}\frac{(n/m)}{2N}, -|\Delta|\frac{(n/m)^2}{4N}\bigg), \textup{ for } n > 0,\\
(-4\pi n)^{-\frac{k}{2}}\sum\limits_{m\mid n}\frac{\begin{pmatrix}\frac{\Delta}{m}\end{pmatrix}}{m\big(4\pi\big|\Delta\frac{(n/m)^2}{4N}\big|\big)^{\frac{k}{4}}}c\bigg(\mathfrak{r}\frac{(n/m)}{2N}, -|\Delta|\frac{(n/m)^2}{4N}\bigg), \textup{ for } n < 0.
\end{cases}
\end{align*}
\end{corollary}
\begin{proof}
Recall that for $\rho_L$ the Weil representation for $L \subset (V,q)$ as in Section \ref{sig21}, the dual Weil representation $\overline{\rho}_L$ is given by that for $L^{-}$, where $L^{-}$ is the lattice $L$ considered as a subset of $(V,-q)$, a rational vector space of signature $(1,2)$. Note furthermore that any function $F(\tau) \in U_{\frac{k+1}{2}, \overline{\tilde{\rho}}_L}(s)$ has a Fourier series of the form
\begin{align}\label{bobs}
\sum\limits_{h \in L'/L}\big[\sum\limits_{\substack {n \in \mathbb{Z} - \textup{sign}(\Delta)q(h) \\n \neq 0}}c(h, n)(4\pi |n|v)^{-\frac{k+1}{4}}W_{\textup{sign}(n)\frac{k+1}{4}, s- \frac{1}{2}}(4\pi |n|v)e(nu)\big]\mathfrak{e}_h.
\end{align}
To consider the Millson-Shimura (Shimura) lift of a Maass cusp form $F(\tau) \in U_{\frac{k+1}{2}, \overline{\tilde{\rho}}_L}(s)$ with $0 \leq k \in 4\mathbb{Z}+2$ ($k \in 4\mathbb{Z}$), one has to regard 
\begin{align*}
\Psi_{\Delta}^*(F, z) = R_{0,z}^{k/2}\int\limits_{\mathcal{F}}\langle L^{k/4}_{\frac{k+1}{2}, \tau}F(\tau) , \Theta^*_{M, \Delta}(z, \tau) \rangle v^{-1/2}d\mu(\tau),
\end{align*}
(respectively $\Phi_{\Delta}^*(F, z)$ and the integral against $\Theta^*_{S, \Delta}(z, \tau)$) with the differential operators adapted accordingly when $k < 0$. By the same arguments as in Proposition \ref{ka}, this is a Maass cusp form of weight $k$. To find the Fourier expansion of $\Psi^*_{\Delta}(F, z) (\Phi^*_{\Delta}(F, z))$, we plug switched theta function into the proof of Theorem \ref{MShim}. The rest follows in the same way as Theorem \ref{MShim}.
\end{proof}

\subsection{Hecke Properties of the Shimura and Millson-Shimura Lift}
We show that the Shimura and Millson-Shimura lift map Hecke eigenforms to Hecke eigenforms.
\begin{theorem} \label{eigenform}
For $k \in 4 \mathbb{Z}$ and an eigenform $F(\tau) \in U_{\frac{k+1}{2}, \tilde{\rho}_L}(s)$ $($or $F(\tau) \in U_{\frac{k+1}{2}, \overline{\tilde{\rho}}_L}(s)$$)$ for all $T_{p^2}$ where $(p,2N) = 1$ and Fourier expansion as in Definition \ref{liftsec} $($as in \eqref{bobs}$)$, the Shimura lift $\Phi_{\Delta}(F, z)$ $($the Millson-Shimura lift $\Psi^*_{\Delta}(F, z)$$)$ is an eigenform for all $T_n$ with $(n, 2N) = 1$.\\ %under the same eigenvalue.
For $k \in 4 \mathbb{Z}+2$ and for an eigenform $F(\tau) \in U_{\frac{k+1}{2}, \tilde{\rho}_L}(s)$ $($$U_{\frac{k+1}{2}, \overline{\tilde{\rho}}_L}(s)$$)$ for all $T_{p^2}$ where $(p,2N) = 1$ and Fourier expansion as in Definition \ref{liftsec} $($as in \eqref{bobs}$)$, the Millson-Shimura lift $\Psi_{\Delta}(F, z)$ $($the Shimura lift $\Phi^*_{\Delta}(F, z)$$)$ is a Hecke eigenform for all $T_n$ with $(n, 2N) = 1$. %under the same eigenvalue.
\end{theorem}
\begin{proof}
We only show this for the Millson-Shimura lift, as the other cases follow in the same way as here or in \cite{MR1244668}. We start with $F(\tau)$ being of weight $-1/2$. We denote the Fourier coefficients of $F(\tau)$ by $c(h, n)$, the $a(n)$ as in Theorem \ref{MShim} and by $u(p)$ the Hecke eigenvalue corresponding to $T_{p^2}$. Assume that $(p,2N) = 1$. We want to prove
\begin{align}\label{tp}
pu(p) a(n) = a(pn) + p^{-1}a(n/p),
\end{align}
We consider the two cases:\\
\textbf{Case 1}: Proving \eqref{tp} for all $a(n)$ with $(n,p) =1$. We have
\begin{align*}
a(n) = \sum\limits_{m\mid n}m^{-2}n\begin{pmatrix}\frac{\Delta}{m}\end{pmatrix}c\bigg(\mathfrak{r}\frac{(n/m)}{2N}, |\Delta|\frac{(n/m)^2}{4N}\bigg)
\end{align*}
and
\begin{align*}
c(ph, np^2) = u(p)c(h,n) -p^{-2}\begin{pmatrix}\frac{\textup{sign}(\Delta)4Nn}{p}\end{pmatrix}c(n,h) -p^{-3}c(h/p, n/p^2).
\end{align*}
Thus, we obtain
\begin{align*}
&a(pn) = \sum\limits_{mk = np}m^{-2}np\begin{pmatrix}\frac{\Delta}{m}\end{pmatrix}c\bigg(\mathfrak{r}\frac{k}{2N}, |\Delta|\frac{k^2}{4N}\bigg)\\
&=  \sum\limits_{m(k'p) = np}m^{-2}np\begin{pmatrix}\frac{\Delta}{m}\end{pmatrix}c\bigg(\mathfrak{r}\frac{k'p}{2N}, |\Delta|\frac{(k'p)^2}{4N}\bigg) \\
&+ \sum\limits_{(m'p)k = np}(mp)^{-2}np\begin{pmatrix}\frac{\Delta}{mp}\end{pmatrix}c\bigg(\mathfrak{r}\frac{k}{2N}, |\Delta|\frac{k^2}{4N}\bigg)\\
&= \sum\limits_{m'k' = n}m^{-2}n\begin{pmatrix}\frac{\Delta}{m}\end{pmatrix}c\bigg(\mathfrak{r}\frac{k}{2N}, |\Delta|\frac{k^2}{4N}\bigg)[pu(p) - p^{-1}\begin{pmatrix}\frac{\Delta}{p}\end{pmatrix}+p^{-1}\begin{pmatrix}\frac{\Delta}{p}\end{pmatrix}]\\
& = pu(p)a(n).
\end{align*}
\textbf{Case 2}: We consider $a(n)$ with $p^{e} \mid n$ and $p^{e+1}\nmid n$. We have
\begin{align*}
a(n) = \sum\limits_{mk = \frac{n}{p^e}}\sum\limits_{j=0}^e c\bigg(\mathfrak{r}\frac{kp^j}{2N}, |\Delta|\frac{k^2p^{2j}}{4N}\bigg) \begin{pmatrix}\frac{\Delta}{mp^{e-j}} \end{pmatrix}\frac{n}{p^e}p^{e}(mp^{e-j})^{-2}.
\end{align*}
We define
\begin{align*}
A_j= c\bigg(\mathfrak{r}\frac{kp^j}{2N}, |\Delta|\frac{k^2p^{2j}}{4N}\bigg),
\end{align*}
which gives us
\begin{align*}
A_1 = (u(p) -p^{-2}\begin{pmatrix}\frac{\Delta}{p} \end{pmatrix})A_0
\textup{ and } A_j = u(p) A_{j-1}-p^{-3}A_{j-2}.
\end{align*}
We obtain $A_j = P_j(u(p))A_0$, where we set
\begin{align*}
P_0(x) = 1, \hspace{4mm} P_1(x) = x -p^{-2}\begin{pmatrix}\frac{\Delta}{p} \end{pmatrix}, \hspace{4mm} P_2(x) = x^2 -p^{-2}\begin{pmatrix}\frac{\Delta}{p} \end{pmatrix}x -p^{-3},
\end{align*}
and
\begin{align*}
P_{j+1}(x)+p^{-3}P_{j-1}(x)=xP_j(x).
\end{align*}
From this, we obtain
\begin{align*}
a(n) = a\Big(\frac{n}{p^e}\Big)p^{-e} \sum\limits_{mk = \frac{n}{p^e}}\sum\limits_{j=0}^e P_j(u(p))\begin{pmatrix}\frac{\Delta}{p^{e-j}} \end{pmatrix}p^{2j}.
\end{align*}
Similarly, we also obtain
\begin{align*}
a\Big(\frac{n}{p}\Big)= a\Big(\frac{n}{p^e}\Big)p^{-e+1}  \sum\limits_{mk = \frac{n}{p^e}}\sum\limits_{j=0}^{e-1} P_j(u(p))\begin{pmatrix}\frac{\Delta}{p^{e-1-j}} \end{pmatrix}p^{2j},
\end{align*}
and
\begin{align*}
a(pn) =  a\Big(\frac{n}{p^e}\Big)p^{-e-1}  \sum\limits_{mk = \frac{n}{p^e}}\sum\limits_{j=0}^{e+1} P_j(u(p))\begin{pmatrix}\frac{\Delta}{p^{e+1-j}} \end{pmatrix}p^{2j}.
\end{align*}
An index shift as in the proof of Proposition 4.1 of \cite{MR1244668}, then gives
\begin{align*}
&a(pn) + p^{-1}a\Big(\frac{n}{p}\Big)\\
&= a\Big(\frac{n}{p^{e}}\Big)p^{-e} \sum\limits_{mk = \frac{n}{p^e}}\sum\limits_{j=1}^e\Big[ p^{-2}P_{j-1}(u(p))+pP_{j+1}\big]\begin{pmatrix}\frac{\Delta}{p^{e-j}} \end{pmatrix}p^{2j}\\
&+ \begin{pmatrix}\frac{\Delta}{p^{e+1}} \end{pmatrix}p^{-1}P_0 +  \begin{pmatrix}\frac{\Delta}{p^{e}} \end{pmatrix}pP_1(u(p))\\
&= pu(p) a(n).
\end{align*}
Due to Hecke eigenfunctions mapping to Hecke eigenfunctions under $L_k$ and $R_k$ (see \cite{str}), this concludes the proof.
\end{proof}

\section{The Vector Valued Katok-Sarnak Formula}\label{kks}
In this section, we state and prove the main theorem of this paper. We will split our results into two different theorems, providing formulas for two slightly different cases, since this form makes it easier to read and comprehend our statement and the proof of it.
The split of our Katok-Sarnak formula differentiates between the cases $Shim$ and $Mish$ as given in the following. The significance of this particular split is that the formula for one of these cases (see Theorem \ref{KSmain}) implies the scalar valued results, see \cite{MR3549627} and \cite{MR4389022}.\\
We denote the space of Maass cusp forms of weight $k$, level $N$ and spectral point $s$ by $U_{k,N}(s)$ and the space of vector valued Maass cusp forms of weight $k$ and spectral point $s$ for $\textup{Mp}_2(\mathbb{Z})$ and the Weil representation $\rho_L$ associated to a lattice $L$ of level $N$ by $U_{k,\rho_L}(s)$. For any $f(z) \in U_{k,N}(s)$ with Fourier coefficients given by $a(n)$ with $a(1) \neq 0$, we define
\begin{align*}
f^{no}(z) := \frac{1}{a(1)}f(z).
\end{align*}
For $F \in U_{k,\rho_L}(s)$ or $F \in U_{k,\overline{\rho}_L}(s)$ and $\Delta \equiv \mathfrak{r}^2 (4N)$ a fundamental discriminant, we set (in case this is well defined)
\begin{itemize}
\item $Shim(F)$ as the normalized Shimura lift for $\Delta > 0$
\begin{align*}
Shim(F) := \begin{cases}
\Phi_{\Delta}^{no}(F, z) &\textup{ for } k \in 4\mathbb{Z} \textup{ and } \Delta > 0,\\
\Phi^{*,no}_{\Delta}(F, z) &\textup{ for } k\in 4\mathbb{Z}+2 \textup{ and } \Delta > 0,
\end{cases}
\end{align*}
and as the normalized Millson-Shimura lift for $\Delta < 0$
\begin{align*}
Shim(F) = \begin{cases}
\Psi_{\Delta}^{no}(F, z) &\textup{ for } k \in 4\mathbb{Z}+2 \textup{ and } \Delta < 0,\\
\Psi^{*,no}_{\Delta}(F, z) &\textup{ for } k\in 4\mathbb{Z} \textup{ and } \Delta < 0.
\end{cases}
\end{align*}
\item $Mish(F)$ as the normalized Millson-Shimura lift for $\Delta > 0$
\begin{align*}
Mish(F) = \begin{cases}
\Psi^{no}_{\Delta}(F, z) &\textup{ for } k \in 4\mathbb{Z}+2 \textup{ and } \Delta > 0,\\
\Psi^{*,no}_{\Delta}(F, z) &\textup{ for } k\in 4\mathbb{Z} \textup{ and } \Delta > 0,
\end{cases}
\end{align*}
and as the normalized Shimura lift for $\Delta < 0$
\begin{align*}
 Mish(F) := \begin{cases}
\Phi_{\Delta}^{no}(F, z) &\textup{ for } k \in 4\mathbb{Z} \textup{ and } \Delta < 0,\\
\Phi^{*,no}_{\Delta}(F, z) &\textup{ for } k\in 4\mathbb{Z}+2 \textup{ and } \Delta <0.
\end{cases}
\end{align*}
\end{itemize}
This is well defined, in case the first Fourier coefficient of  $Shim$ (or $Mish$) does not equal to $0$. In particular, this is the case, when it is a non-zero multiple of a newform.
Note that we refrain from indicating the dependency of $Shim$ and $Mish$ on $\Delta$ and $k$ by an index, due to the simpler notation. However, it would technically be more correct to denote them by $Shim_{k,\Delta}$ and $Mish_{k,\Delta}$.

\begin{theorem}[The Katok-Sarnak Formula I]\label{KSmain}

Let $k \in 4 \mathbb{Z}$ and $\varphi \in U_{k,N}(s)$ be an even normalized newform. Denote by $\{F_j\}$ an orthonormal basis of $U_{\frac{k+1}{2}, \rho_L}\left(\frac{s}{2}+\frac{1}{4}\right)$ which consists of simultaneous eigenforms for all Hecke operators $T_{p^2}$ for $(p,2N) = 1$ and by $b_j(h, n)$ the $n$-th Fourier coefficient corresponding to the $h$-th component of $F_j$. Let $\mathfrak{r}\in\mathbb{Z}$, $\Delta \equiv \mathfrak{r}^2 (4N)$ be a fundamental discriminant and $D \in \mathbb{Q}$ such that $4ND$ is a square modulo $4N$. For $k\geq 0$, we have
\begin{align*}
&\sum\limits_{Shim(F_j) =\varphi} (-1)^{k/4}2|\Delta D|^{\frac{1}{2}}(4N)^{\frac{k}{4}} \overline{b_j(\mathfrak{r}/2N, \Delta/4N)}b_j(h, D)\\
&= (\varphi, \varphi)^{-1}|\Delta D|^{\frac{k}{4}}\begin{cases}
\frac{\sqrt{\pi}}{2} Tr^{\textup{s}(\Delta)}_{\Delta, \mathfrak{r}}\left(L_{k,z}^{k/2}\varphi, h, D \right)& \textup{ for } \Delta D < 0,\\
2^{\frac{k}{2}-1}Tr_{\Delta, \mathfrak{r}, k}\left(\varphi, h, D\right) & \textup{ for } \Delta,D > 0,\\
 \frac{1}{4}2^{-\frac{k}{2}}\frac{\Gamma(s+k/2)}{\Gamma(s-k/2)}\cdot Tr_{\Delta, \mathfrak{r},k}(L_{k,z}\varphi, h, -D) & \textup{ for } \Delta, D < 0,
\end{cases}
\end{align*}
where for any weight $0$ Maass cusp form $g$ for $\Gamma_0(N)$, we set
\begin{align*}
&Tr^{\textup{s}(\Delta)}_{\Delta, \mathfrak{r}}\left(g, h, D \right)\\
&=\left[Tr^+_{\Delta, \mathfrak{r}}\left(g, h, \textup{sign}(\Delta)D \right) +(-1)^{\frac{1}{2}- \textup{sign}(\Delta)\frac{1}{2}} Tr^-_{\Delta, \mathfrak{r}}\left(g, h, \textup{sign}(\Delta)D \right)\right].
\end{align*}
For $k < 0$, we have
\begin{align*}
&\sum\limits_{Shim(F_j) =\varphi} (-1)^{k/4}2|\Delta D|^{\frac{1}{2}}(4N)^{\frac{k}{4}} \overline{b_j(\mathfrak{r}/2N, \Delta/4N)}b_j(h, D)\\
&= (\varphi, \varphi)^{-1}|\Delta D|^{\frac{k}{4}} \begin{cases}
\frac{\sqrt{\pi}}{2}\frac{\Gamma\left(s+k/2\right)}{\Gamma\left(s-k/2\right)}Tr^{\textup{s}(\Delta)}_{\Delta, \mathfrak{r}}\left(R_{k,z}^{-k/2}\varphi, h, D \right) & \textup{ for } \Delta D < 0,\\
2^{\frac{k}{2}-1}Tr_{\Delta, \mathfrak{r}, k}\left(\varphi, h, D\right) & \textup{ for } \Delta,D > 0,\\
 \frac{1}{4}2^{-\frac{k}{2}}\frac{\Gamma(s+k/2)}{\Gamma(s-k/2)}\cdot Tr_{\Delta, \mathfrak{r},k}(R_{k,z}\varphi, h, -D) & \textup{ for } \Delta, D < 0.
\end{cases}
\end{align*}
Let $k \in 4 \mathbb{Z}+2$ and $\varphi \in U_{k,N}(s)$ be an even normalized newform. Denote by $\{F_j\}$ an orthonormal basis of $U_{\frac{k+1}{2}, \overline{\rho}_L}\left(\frac{s}{2}+\frac{1}{4}\right)$ which consists of simultaneous eigenforms for all $T_{p^2}$ for $(p,2N) = 1$ and by $b_j(h, n)$ the $n$-th Fourier coefficient corresponding to the $h$-th component of $F_j$. Let $\mathfrak{r}\in\mathbb{Z}$, $\Delta \equiv \mathfrak{r}^2 (4N)$ be a fundamental discriminant and $D \in \mathbb{Q}$ such that $-4ND$ is a square modulo $4N$. For $k>0$, we have
\begin{align*}
&\sum\limits_{Shim(F_j) = \varphi}(-1)^{\frac{k+2}{4}}2|\Delta D|^{\frac{1}{2}}(4N)^{\frac{k}{4}} \overline{b_j(\mathfrak{r}/2N, -\Delta/4N)}b_j(h, D)\\
&= (\varphi, \varphi)^{-1}|\Delta D|^{\frac{k}{4}}\begin{cases}
\frac{\sqrt{\pi}}{2}Tr^{\textup{s}(\Delta)}_{\Delta, \mathfrak{r}}\left(L_{k,z}^{k/2}\varphi, h, D \right) & \textup{ for } \Delta D < 0,\\
2^{\frac{k}{2}-1}Tr_{\Delta, \mathfrak{r}, k}\left(\varphi, h, D\right) & \textup{ for } \Delta,D > 0,\\
 \frac{1}{4}2^{-\frac{k}{2}}\frac{\Gamma(s+k/2)}{\Gamma(s-k/2)}\cdot Tr_{\Delta, \mathfrak{r},k}(L_{k,z}\varphi, h, -D) & \textup{ for } \Delta, D < 0,
\end{cases}
\end{align*}
And for $k < 0$, we have
\begin{align*}
&\sum\limits_{Shim(F_j) =\varphi} (-1)^{\frac{k+2}{4}}2|\Delta D|^{\frac{1}{2}}(4N)^{\frac{k}{4}} \overline{b_j(\mathfrak{r}/2N, -\Delta/4N)}b_j(h, D)\\
&= (\varphi, \varphi)^{-1}|\Delta D|^{\frac{k}{4}} \begin{cases}
\frac{\sqrt{\pi}}{2}\frac{\Gamma\left(s+k/2\right)}{\Gamma\left(s-k/2\right)}Tr^{\textup{s}(\Delta)}_{\Delta, \mathfrak{r}}\left(R_{k,z}^{-k/2}\varphi, h, D \right) & \textup{ for } \Delta D < 0,\\
2^{\frac{k}{2}-1}Tr_{\Delta, \mathfrak{r}, k}\left(\varphi, h, D\right) & \textup{ for } \Delta,D > 0,\\
 \frac{1}{4}2^{-\frac{k}{2}}\frac{\Gamma(s+k/2)}{\Gamma(s-k/2)}\cdot Tr_{\Delta, \mathfrak{r},k}(R_{k,z}\varphi, h, -D) & \textup{ for } \Delta, D < 0.
\end{cases}
\end{align*}
\end{theorem}

\begin{proof}[Proof of Theorem 1.1]
To each even normalized newform $\varphi \in U_{k,N}(s)$, we can associate a Maass cusp form $I_{\Delta}(\tau, \varphi)$ of weight $\frac{k+1}{2}$, given by
\begin{align*}
I_{\Delta}(\tau, \varphi)
:= \begin{cases}
I_{S,\Delta}(\tau, \varphi) & \textup{ if } \Delta > 0 \textup{ and } k \in 4\mathbb{Z},\\
I^*_{S,\Delta}(\tau, \varphi) & \textup{ if } \Delta > 0 \textup{ and } k \in 4\mathbb{Z}+2,\\
I^*_{M,\Delta}(\tau, \varphi) & \textup{ if } \Delta < 0 \textup{ and } k \in 4\mathbb{Z},\\
I_{M,\Delta}(\tau, \varphi) & \textup{ if } \Delta < 0 \textup{ and } k \in 4\mathbb{Z}+2.\\
\end{cases}
\end{align*}
Denote the Fourier series of $I_{\Delta}(\tau,\varphi)$ by
\begin{align*}
&I_{\Delta}(\tau,\varphi)\\
&= \sum\limits_{h \in L'/L}\big[\sum\limits_{\substack {D \in \mathbb{Z}+\tilde{q}(h)\\D \neq 0}}\tilde{C}(h, D)(4\pi |D|v)^{-\frac{k+1}{4}}W_{\textup{sign}(D)\frac{k+1}{4}, \frac{s}{2}- \frac{1}{4}}(4\pi |D|v)e(Du)\big]\mathfrak{e}_h,
\end{align*}
where $\tilde{q} = (-1)^{k/2}q$. There are now two possible methods to compute the Fourier coefficients $\tilde{C}(h,D)$.
The first method is described in Section \ref{Siegelcoeff}. The $\tilde{C}(h,D)$ are then given by Theorem \ref{fünf}, Theorem \ref{sechs}, Theorem \ref{sieben} and Theorem \ref{acht}. \\
The second method works as follows: Denote
\begin{align*}
\rho := \begin{cases}
\rho_L & \textup{ for } k \in 4\mathbb{Z},\\
\overline{\rho}_L & \textup{ for } k \in 4\mathbb{Z}+2.
\end{cases}
\end{align*}
Pick then a basis $\{F_j\}$ of $U_{\frac{k+1}{2}, \rho}(\frac{s}{2}+\frac{1}{4})$ which consists of simultaneous eigenforms for all Hecke operators $T_{p^2}$ with $(n,2N) =1$ . This exists due to \cite{MR2574974}, Theorem 4.12. Thus, we can represent $I_{\Delta}(\tau, \varphi)$ in the following way
\begin{align*}
I_{\Delta}(\tau, \varphi) = \sum\limits_{j}(I_{\Delta} , F_j ) F_j(\tau),
\end{align*}
where we write $(\cdot, \cdot)$ for the Petersson inner product on $U_{\frac{k+1}{2}, \rho}(s)$. Denote by $b_j(h, n)$ the $n$-th Fourier coefficient of the $h$-th component of $F_j$. The $D$-th Fourier coefficient of the $h$-th component of $I_{\Delta}(\tau, \varphi)$ is then given by
\begin{align*}
\tilde{C}(h,D) = \sum\limits_{j}( I_{\Delta}, F_j ) b_j(h, D).
\end{align*}
For the next step, we only consider the case that $0<k \in 4\mathbb{Z}$ and $\Delta > 0$, since any other case follows by the exact same arguments and one simply has to adapt the theta function and the application of the lowering and raising operator. It remains to compute
\begin{align*}
&( I_{\Delta}, F_j ) = \int\limits_{\mathcal{F}}\langle I_{\Delta}(\tau, \varphi),F_j(\tau)\rangle v^{(k+1)/2}d\mu(\tau)\\
&= \int\limits_{\mathcal{F}}\bigg\langle\bigg[R_{1/2, \tau}^{k/4}\int\limits_{\Gamma\setminus\mathbb{H}}L_{k,z}^{k/2}\varphi(z) \Theta_{S, \Delta}(z,\tau)d\mu(z)\bigg], F_j(\tau)\bigg\rangle v^{(k+1)/2}d\mu(\tau).
\end{align*}
By a Stokes argument, see Lemma 4.2 of \cite{MR1903920} or Lemma 3.11 of \cite{hoe}, we obtain
\begin{align*}
& \int\limits_{\mathcal{F}}\bigg\langle\bigg[R_{1/2, \tau}^{k/4}\int\limits_{\Gamma\setminus\mathbb{H}}L_{k,z}^{k/2}\varphi(z) \Theta_{S, \Delta}(z, \tau)d\mu(z)\bigg], F_j(\tau)\bigg\rangle v^{(k+1)/2}d\mu(\tau)\\
&= (-1)^{k/4} \int\limits_{\mathcal{F}}\bigg\langle\bigg[\int\limits_{\Gamma\setminus\mathbb{H}}L_{k,z}^{k/2}\varphi(z) \Theta_{S, \Delta}(z, \tau)d\mu(z)\bigg], L_{\frac{k+1}{2}, \tau}^{k/4}F_j(\tau)\bigg\rangle v^{1/2}d\mu(\tau).
\end{align*}
By a Fubini argument, we change the order of integration, which gives us that the term above is equal to
\begin{align*}
(-1)^{k/4} \int\limits_{\Gamma\setminus\mathbb{H}}L_{k,z}^{k/2}\varphi(z)\overline{\bigg[\int\limits_{\mathcal{F}} \bigg\langle L_{\frac{k+1}{2}, \tau}^{k/4}F_j(\tau), \Theta_{S, \Delta}(z, \tau)\bigg\rangle v^{1/2}d\mu(\tau)\bigg]}d\mu(z),
\end{align*}
which, again by a Stokes argument and the fact that $4\mid k$, is equal to
\begin{align*}
&(-1)^{k/4} \int\limits_{\Gamma\setminus\mathbb{H}}\varphi(z)\overline{R_{0,z}^{k/2}\bigg[\int\limits_{\mathcal{F}} \bigg\langle L_{\frac{k+1}{2}, \tau}^{k/4}F_j(\tau), \Theta_{S,\Delta}(z, \tau),\bigg\rangle v^{1/2}d\mu(\tau)\bigg]}y^kd\mu(z)\\
&= (-1)^{k/4}\big(\varphi, \Phi_{\Delta}(F_j, z)\big),
\end{align*}
where $\Phi_{\Delta}(F_j, z)$ is the Shimura lift of $F_j$. \\
Executing this step for all cases of $I_{\Delta}(\tau, \varphi)$ (i.e.\ all $\Delta$ and $k$), we obtain that $( I_{\Delta}, F_j ) = \big(\varphi, \phi_{\Delta}(F_j, z)\big)$ with
\begin{align*}
\phi_{\Delta}(F_j, z) := \begin{cases}
\Phi_{\Delta}(\tau, \varphi) & \textup{ if } \Delta > 0 \textup{ and } k \in 4\mathbb{Z},\\
\Phi^*_{\Delta}(\tau, \varphi) & \textup{ if } \Delta > 0 \textup{ and } k \in 4\mathbb{Z}+2,\\
\Psi^*_{M,\Delta}(\tau, \varphi) & \textup{ if } \Delta < 0 \textup{ and } k \in 4\mathbb{Z},\\
\Psi_{M,\Delta}(\tau, \varphi) & \textup{ if } \Delta < 0 \textup{ and } k \in 4\mathbb{Z}+2.\\
\end{cases}
\end{align*}
We know that $\phi_{\Delta}(F_j, z)$ is a weight $k$ even eigenform for all Hecke operators $T_n$ with $(n, 2N) = 1$ due to Theorem \ref{eigenform}. It follows by a multiplicity one argument that $\phi_{\Delta}(F_j, z)$ is the multiple of an even newform of level $N$, which we call $\tilde{\varphi}_j(z)$ . In order for $\big(\varphi, \phi_{\Delta}(F_j, z)\big)$ to not automatically be zero, then, the function $\tilde{\varphi}_j(z)$ must be multiple of $\varphi(z)$. Hence, we have
\begin{align*}
\big(\varphi, \phi_{\Delta}(F_j, z)\big) =\begin{cases} \overline{c}\big(\varphi(z), \varphi(z)\big) &\textup{ if }\phi_{\Delta}(F_j, z) = c\varphi(z),\\
0 & \textup{ otherwise },
\end{cases}
\end{align*}
where, since $\varphi$ is normalized, the constant $c$ is given by the first Fourier coefficient of $\phi_{\Delta}(F_j, z)$ (which is non-zero, due to $\phi_{\Delta}(F_j, z)$ being a newform). This is given by Theorem \ref{MShim} and Corollary \ref{cs}. Comparing both results for the $\tilde{C}(h, D)$ then gives us the formula.
\end{proof}

There is a second Katok-Sarnak formula that is associated to $Mish$, which we did not give in the introduction. It does not correspond to the classic scalar valued formulas of \cite{MR1244668}, \cite{MR3549627} and \cite{MR4389022}.

\begin{theorem}[The Katok-Sarnak Formula II]\label{KSsub}

Let $k \in 4 \mathbb{Z}$ and $\varphi \in U_{k,N}(s)$ be an even normalized newform. Denote by $\{F_j\}$ an orthonormal basis of $U_{\frac{k+1}{2}, \overline{\rho}_L}\left(\frac{s}{2}+\frac{1}{4}\right)$ which consists of simultaneous eigenforms for all Hecke operators $T_{p^2}$ for $(p,2N) = 1$ and by $b_j(h, n)$ the $n$-th Fourier coefficient corresponding to the $h$-th component of $F_j$. Let $\mathfrak{r} \in \mathbb{Z}$, $\Delta \equiv \mathfrak{r}^2 (4N)$ be a fundamental discriminant and $D \in \mathbb{Q}$ such that $-4ND$ is a square modulo $4N$. For $k\geq0$, we have
\begin{align*}
&\sum\limits_{Mish(F_j) =\varphi} (-1)^{k/4}2(-i)|\Delta D|^{\frac{1}{2}}(4N)^{\frac{k}{4}} \overline{b_j(\mathfrak{r}/2N, -\Delta/4N)}b_j(h, D)\\
&= (\varphi, \varphi)^{-1}|\Delta D|^{\frac{k}{4}}\begin{cases}
\frac{\sqrt{\pi}}{2} Tr^{\textup{ß}(\Delta)}_{\Delta, \mathfrak{r}}\left(L_{k,z}^{k/2}\varphi, h, D \right)& \textup{ for } \Delta D < 0,\\
2^{\frac{k}{2}-1}Tr_{\Delta, \mathfrak{r}, k}\left(\varphi, h, D\right) & \textup{ for } \Delta,D > 0,\\
 -\frac{1}{4}2^{-\frac{k}{2}}\frac{\Gamma(s+k/2)}{\Gamma(s-k/2)}\cdot Tr_{\Delta, \mathfrak{r},k}(L_{k,z}\varphi, h, -D) & \textup{ for } \Delta, D < 0,
\end{cases}
\end{align*}
where for any weight $0$ Maass cusp form $g$ for $\Gamma_0(N)$, we set
\begin{align*}
&Tr^{\textup{ß}(\Delta)}_{\Delta, \mathfrak{r}}\left(g, h, D \right)\\
&=\left[Tr^-_{\Delta, \mathfrak{r}}\left(g, h, -\textup{sign}(\Delta)D \right) +(-1)^{\frac{1}{2}+ \textup{sign}(\Delta)\frac{1}{2}} Tr^+_{\Delta, \mathfrak{r}}\left(g, h, -\textup{sign}(\Delta)D \right)\right].
\end{align*}
For $k < 0$, we have
\begin{align*}
&\sum\limits_{Mish(F_j) =\varphi} (-1)^{k/4}2(-i)|\Delta D|^{\frac{1}{2}}(4N)^{\frac{k}{4}} \overline{b_j(\mathfrak{r}/2N, -\Delta/4N)}b_j(h, D)\\
&= (\varphi, \varphi)^{-1}|\Delta D|^{\frac{k}{4}} \begin{cases}
\frac{\sqrt{\pi}}{2}\frac{\Gamma\left(s+k/2\right)}{\Gamma\left(s-k/2\right)}Tr^{\textup{ß}(\Delta)}_{\Delta, \mathfrak{r}}\left(R_{k,z}^{-k/2}\varphi, h, D \right) & \textup{ for } \Delta D < 0,\\
2^{\frac{k}{2}-1}Tr_{\Delta, \mathfrak{r}, k}\left(\varphi, h, D\right) & \textup{ for } \Delta,D > 0,\\
- \frac{1}{4}2^{-\frac{k}{2}}\frac{\Gamma(s+k/2)}{\Gamma(s-k/2)}\cdot Tr_{\Delta, \mathfrak{r},k}(R_{k,z}\varphi, h, -D) & \textup{ for } \Delta, D < 0.
\end{cases}
\end{align*}
Let $k \in 4 \mathbb{Z}+2$ and $\varphi \in U_{k,N}(s)$ be an even normalized newform. Denote by $\{F_j\}$ an orthonormal basis of $U_{\frac{k+1}{2}, \rho_L}\left(\frac{s}{2}+\frac{1}{4}\right)$ which consists of simultaneous eigenforms for all $T_{p^2}$ for $(p,2N) = 1$ and by $b_j(h, n)$ the $n$-th Fourier coefficient corresponding to the $h$-th component of $F_j$. Let $\mathfrak{r} \in \mathbb{Z}$, $\Delta \equiv \mathfrak{r}^2 (4N)$ be a fundamental discriminant and $D \in \mathbb{Q}$ such that $4ND$ is a square modulo $4N$. For $k>0$, we have
\begin{align*}
&\sum\limits_{Mish(F_j) = \varphi}(-1)^{\frac{k+2}{4}}2(-i)|\Delta D|^{\frac{1}{2}}(4N)^{\frac{k}{4}} \overline{b_j(\mathfrak{r}/2N, \Delta/4N)}b_j(h, D)\\
&= (\varphi, \varphi)^{-1}|\Delta D|^{\frac{k}{4}}\begin{cases}
\frac{\sqrt{\pi}}{2}Tr^{\textup{ß}(\Delta)}_{\Delta, \mathfrak{r}}\left(L_{k,z}^{k/2}\varphi, h, D \right) & \textup{ for } \Delta D < 0,\\
2^{\frac{k}{2}-1}Tr_{\Delta, \mathfrak{r}, k}\left(\varphi, h, D\right) & \textup{ for } \Delta,D > 0,\\
 -\frac{1}{4}2^{-\frac{k}{2}}\frac{\Gamma(s+k/2)}{\Gamma(s-k/2)}\cdot Tr_{\Delta, \mathfrak{r},k}(L_{k,z}\varphi, h, -D) & \textup{ for } \Delta, D < 0.
\end{cases}
\end{align*}
And for $k < 0$, we have
\begin{align*}
&\sum\limits_{Mish(F_j) = \varphi}(-1)^{\frac{k+2}{4}}2(-i)|\Delta D|^{\frac{1}{2}}(4N)^{\frac{k}{4}} \overline{b_j(\mathfrak{r}/2N, \Delta/4N)}b_j(h, D)\\
&= (\varphi, \varphi)^{-1}|\Delta D|^{\frac{k}{4}} \begin{cases}
\frac{\sqrt{\pi}}{2}\frac{\Gamma\left(s+k/2\right)}{\Gamma\left(s-k/2\right)}Tr^{\textup{ß}(\Delta)}_{\Delta, \mathfrak{r}}\left(R_{k,z}^{-k/2}\varphi, h, D \right) & \textup{ for } \Delta D < 0,\\
2^{\frac{k}{2}-1}Tr_{\Delta, \mathfrak{r}, k}\left(\varphi, h, D\right) & \textup{ for } \Delta,D > 0,\\
 -\frac{1}{4}2^{-\frac{k}{2}}\frac{\Gamma(s+k/2)}{\Gamma(s-k/2)}\cdot Tr_{\Delta, \mathfrak{r},k}(R_{k,z}\varphi, h, -D) & \textup{ for } \Delta, D < 0.
\end{cases}
\end{align*}
\end{theorem}
\begin{proof}
This time, to each $\varphi \in U_{k,N}(s)$ we associate a Maass cusp form $J_{\Delta}(\tau, \varphi)$ of weight $\frac{k+1}{2}$, given by
\begin{align*}
J_{\Delta}(\tau, \varphi)
:= \begin{cases}
I_{S,\Delta}(\tau, \varphi) & \textup{ if } \Delta < 0 \textup{ and } k \in 4\mathbb{Z},\\
I^*_{S,\Delta}(\tau, \varphi) & \textup{ if } \Delta < 0 \textup{ and } k \in 4\mathbb{Z}+2,\\
I^*_{M,\Delta}(\tau, \varphi) & \textup{ if } \Delta > 0 \textup{ and } k \in 4\mathbb{Z},\\
I_{M,\Delta}(\tau, \varphi) & \textup{ if } \Delta > 0 \textup{ and } k \in 4\mathbb{Z}+2.\\
\end{cases}
\end{align*}
We proceed to run through the steps of the proof of Theorem \ref{KSmain} with $J_{\Delta}(\tau, \varphi)$ instead of $I_{\Delta}(\tau, \varphi)$ and $\overline{\rho}$ instead of $\rho$. We end up with having to consider $(\varphi, \psi_{\Delta}(F_j,z))$ instead of $(\varphi, \phi_{\Delta}(F_j,z))$, where
\begin{align*}
\psi_{\Delta}(F_j, z) := \begin{cases}
\Phi_{\Delta}(\tau, \varphi) & \textup{ if } \Delta <0  \textup{ and } k \in 4\mathbb{Z},\\
\Phi^*_{\Delta}(\tau, \varphi) & \textup{ if } \Delta < 0 \textup{ and } k \in 4\mathbb{Z}+2,\\
\Psi^*_{M,\Delta}(\tau, \varphi) & \textup{ if } \Delta > 0 \textup{ and } k \in 4\mathbb{Z},\\
\Psi_{M,\Delta}(\tau, \varphi) & \textup{ if } \Delta > 0 \textup{ and } k \in 4\mathbb{Z}+2.\\
\end{cases}
\end{align*}
Comparing Theorem \ref{fünf}, Theorem \ref{sechs}, Theorem \ref{sieben} and Theorem \ref{acht} with Theorem \ref{MShim} and Corollary \ref{cs} again then concludes the proof.
\end{proof}

\begin{remark}
Our formula, as given in Theorem \ref{KSmain} can be used to derive a scalar valued Katok-Sarnak formulas in \cite{MR3549627} and \cite{MR4389022}. This works by means of an isomorphism given in Theorem 5.6 of \cite{MR781735}.\footnote{Note, that our formula differs from that of \cite{MR4389022} by an additional factor $\frac{\Gamma\left(s+k/2\right)}{\Gamma\left(s-k/2\right)}$ in the case $\Delta D < 0$ for $k <0$ and in the case $\Delta, D < 0$.}
\end{remark}

\bibliographystyle{alpha}
%\chapterstyle{unnumberedchapter}
\bibliography{Kapitel/Appendix/bib.bib}

\end{document}